\documentclass{amsart}

\usepackage{graphicx}

\usepackage{mathtools}

\numberwithin{equation}{section}

\usepackage{amssymb}
\usepackage[UKenglish]{babel}
\usepackage{currfile}
\usepackage[UKenglish]{datetime}
\usepackage{color}
\usepackage{enumerate}
\usepackage{amsthm}
\usepackage{hyperref}

\newcommand{\zeroset}[1]{{\bf Z}(#1)}

\newcommand{\dist}{{\mathrm d}}
\newcommand{\ball}{{\mathrm B}}
\newcommand{\sphere}{\mathrm S}

\newcommand{\cont}{{\mathrm C }}
\newcommand{\contb}{{\mathrm C}_{\mathrm b}}
\newcommand{\contc}{{\mathrm C}_{\mathrm c}}
\newcommand{\conto}{{\mathrm C}_{0}}

\newcommand{\contX}{\cont(X)}
\newcommand{\contbX}{\contb(X)}
\newcommand{\contcX}{\contc(X)}
\newcommand{\contoX}{\conto(X)}

\newcommand{\LipbX}{\mathrm{Lip}_{\mathrm b}(X)}

\newcommand{\orderdual}[1]{{#1}^\sim}
\newcommand{\ordercontn}[1]{{#1}^\sim_{\mathrm n}}
\newcommand{\ordercontc}[1]{{#1}^\sim_{\mathrm c}}

\newcommand{\zerofunction}{{\mathbf 0}}
\newcommand{\onefunction}{{\mathbf 1}}

\newtheorem{thm}{Theorem}[section]
\newtheorem{lem}[thm]{Lemma}
\newtheorem{pro}[thm]{Proposition}

\newtheorem{rem}[thm]{Remark}
\newtheorem{defn}[thm]{Definition}
\newtheorem{cor}[thm]{Corollary}
\newtheorem*{thm*}{Theorem}
\newtheorem*{defn*}{Definition}

\begin{document}
	
\title{On order continuous duals of vector lattices of continuous functions}

\author{Marcel de Jeu}
\address[Marcel de Jeu]{Mathematical Institute, Leiden University, P.O.\ Box 9512, 2300 RA Leiden\\
	the Netherlands\\
	and Department of Mathematics and Applied Mathematics, University of Pretoria, Cor\-ner of Lynnwood Road and Roper Street, Hatfield 0083, Pretoria\\
	South Africa
}
   \email{mdejeu@math.leidenuniv.nl}

\author{Jan Harm van der Walt}
\address[Jan Harm van der Walt]{Department of Mathematics and Applied Mathematics, University of Pretoria, Cor\-ner of Lynnwood Road and Roper Street, Hatfield 0083, Pretoria\\
	South Africa
}
\email{janharm.vanderwalt@up.ac.za}

\begin{abstract}
A topological space $X$ is called resolvable if it contains a dense subset with dense complement.
Using only basic principles, we show that whenever the space $X$ has a resolving subset that can be written as an at most countably infinite union of subsets, in such a way that a given vector lattice of (not necessarily bounded) continuous functions on $X$ separates every point outside the resolving subset from each of its constituents, then the order continuous dual of this lattice is trivial.  In order to apply this result in specific cases, we show that several spaces have resolving subsets that can be written as at most countably infinite unions of closed nowhere dense subsets. An appeal to the main result then yields, for example, that, under appropriate conditions, vector lattices of continuous functions on separable spaces, metric spaces, and topological vector spaces have trivial order continuous duals if they separate points and closed nowhere dense subsets. Our results in this direction extend known results in the literature. We also show that, under reasonably mild separation conditions, vector lattices of continuous functions on locally connected $T_1$ Baire spaces without isolated points have trivial order continuous duals. A discussion of the relation between our results and the non-existence of non-zero normal measures is included.  
\end{abstract}

\keywords{Vector lattice of continuous functions, order continuous dual, resolvable topological space, normal measure}

\subjclass{Primary 46E05; Secondarry 47B65, 4D99}

\maketitle 

\section{Introduction and overview}\label{sec:introduction_and_overview}

Let $E$ be a vector lattice.  We recall that an order bounded linear functional $\varphi$ on $E$ is \emph{order continuous}  if $\inf\{\,\lvert\varphi(f_i)\rvert : i\in I\,\}=0$ whenever $(f_i)_{i\in I}$ is a net in $E$ such that $f_i\downarrow 0$ in $E$.
The vector lattice of all order continuous functionals on $E$ is called the \emph{order continuous dual of $E$}.  This paper is centred around the following result on the triviality of order continuous duals of vector lattices of continuous functions. It will be established in Section~\ref{sec:main_triviality_theorem}.

\begin{thm*}
Let $X$ be a non-empty topological space, and let $E$ be a vector sublattice of the real-valued continuous functions on $X$. Suppose that there exists an at most countably infinite collection $\{\,\Gamma_n : n\in\mathbb{N}\,\}$ of non-empty subsets of $X$ such that
 \begin{enumerate}[(1)]
    \item $\bigcup_{n\in\mathbb{N}}\Gamma_n$ and $X\setminus \bigcup_{n\in\mathbb{N}}\Gamma_n$ are both dense in $X$, and
    \item for every $x\in X\setminus \bigcup_{n\in\mathbb{N}}\Gamma_n$ and every $k\in\mathbb N$, there exists an element $u$ of $E^+$ such that $u(x)=1$ and $u(y)=0$ for all $y\in \Gamma_k$.
\end{enumerate}
Then 0 is the only order continuous linear functional on $E$.

If, in addition, $X$ satisfies the countable chain condition \textup{(}in particular, if, in addition, $X$ is separable\textup{)}, then 0 is also the only $\sigma$-order continuous linear functional on $E$.
\end{thm*}

Let us make a few comments on this theorem.  First of all, we mention explicitly that $E$ need not consist of bounded functions, and that there are no topological properties of $X$ supposed, other than what is in the theorem. Singletons need  not even be closed subsets. In several of the existing results in the literature on the triviality of order continuous duals, $X$ is supposed to be a locally compact Hausdorff space and $E$ is supposed to be a vector lattice of bounded continuous functions on $X$.  More often than not, the Riesz representation theorem is then used to reduce the  study of order continuous  linear functionals on $E$ to the analysis of the structure of the measures that represent them.  Our approach is quite different. Using only first principles,  the argument with nets and sequences is completely elementary, and further properties of $X$ or of the elements of $E$ play no role.

Let us also mention that it is not supposed that the $\Gamma_n$ are pairwise disjoint, or even different. If one wishes, disjointness can nevertheless always be obtained by replacing each $\Gamma_n$ with $\Gamma_n\setminus{\bigcup_{j=1}^{n-1}}\Gamma_j$ for $n\geq 2$, and disregarding all empty sets that might occur in this process.

There are no explicit topological conditions on the $\Gamma_n$ in the theorem, but the coupling with the continuous functions in $E$ implies that some properties are still automatic. For one thing, $\Gamma_k$ is nowhere dense in $X$ for all $k\geq 1$. To see this, fix $k\geq 1$. For every $x\in X\setminus\bigcup_{n\in\mathbb{N}}\Gamma_n$, assumption (2) of the theorem implies that there exists an element $u_x$ of $E$ such that $u_x(x)=1$ and $u_x(y)=0$ for all $y\in\Gamma_k$. Consequently, the open neighbourhood $V_x\coloneqq u^{-1}[(\frac{1}{2},2)]$ of $x$ is contained in $X\setminus\Gamma_k$. Set $V\coloneqq\bigcup_{x\in X\setminus\bigcup_{\in\mathbb{N}}\Gamma_n}V_x$. Then $V$ is an open subset of $X$, and  $X\setminus\bigcup_{n\in\mathbb{N}}\Gamma_n\subseteq V\subseteq X\setminus\Gamma_k$. Since $X\setminus\bigcup_{n\in\mathbb{N}}\Gamma_n$ is dense in $X$, we conclude that $X\setminus\Gamma_k$ contains the open dense subset $V$ of $X$. If $\widetilde V$ were a non-empty open subset of $X$ that is contained in  $\overline{\Gamma}_k$, then $\widetilde V\cap V$ would be a non-empty open subset of $X$ that is contained in $\overline\Gamma_k\setminus\Gamma_k$, which is impossible. Hence $\overline{\Gamma}_k$ has empty interior, i.e.\ $\Gamma_k$ is nowhere dense.

If $X$ is a Baire space, then this implies that one may as well assume that the $\Gamma_n$ are all closed nowhere dense subsets of $X$, simply by replacing $\Gamma_n$ with $\overline{\Gamma}_n$ for all $n\geq 1$. Indeed, $\bigcup_{n\in\mathbb N}\overline\Gamma_n$ is then certainly dense in $X$, but $X\setminus\bigcup_{n\in\mathbb N}\overline\Gamma_n$ is also still dense, since the fact that $X$ is a Baire space now implies that $\bigcup_{n\in\mathbb N}\overline\Gamma_n$ still has empty interior. Furthermore, if $k\geq 1$ is fixed, and if $x\notin\bigcup_{n\in\mathbb N}\overline\Gamma_n$, then $x\notin\bigcup_{n\in\mathbb N}\Gamma_n$, so that there exists an element $u$ of $E^+$ such that $u(x)=1$ and $u(y)=0$ for all $y\in\Gamma_k$. By continuity, we then also have that $u(y)=0$ for all $y\in\overline\Gamma_k$. Hence the $\overline{\Gamma}_n$ also satisfy the assumptions in the theorem. Of course, if the original $\Gamma_n$ are pairwise disjoint, this property may be lost when replacing them with their closures.

In the applications that we give in Section~\ref{sec:triviality_of_order_continuous_duals}, the $\Gamma_n$ are, in fact, typically closed nowhere dense subsets of $X$. The above discussion shows that this is not totally unexpected.

In order to continue our discussion of the theorem, we recall the following notion that was introduced by Hewitt \cite{Hewitt1943}.

\begin{defn}\label{def:resolvable}
A topological space $X$ is called \emph{resolvable} if there exists a subset $D$ of $X$ such that $D$ and $X\setminus D$ are both dense in $X$. We shall then say that such a subset $D$ is a \emph{resolving subset of $X$}, or that it \emph{resolves $X$}.
\end{defn}

Hence part (1) of the hypotheses implies that  $X$ should be resolvable, but not just that: Combination with part (2) shows that it must be resolvable in a special manner.
If one thinks of $E$ as given, then, using that the $\Gamma_n$ may as well be taken to be disjoint, the interpretation of the parts (1) and (2) taken together is the following: There exists a resolving subset of $X$ that can be split into at most countably infinitely many non-empty parts, in such a way that $E$ is sufficiently rich with respect to each of these parts in order to be able to supply all separating functions in part (2).  If one wants to coin a more or less suggestive terminology, then, although it does not capture everything, one could call $\{\,\Gamma_n : n\in\mathbb N\,\}$ an $E$-separated resolution of $X$.

Furthermore, let us observe that, if the hypotheses of the theorem are satisfied for $E$, they are obviously also satisfied for every superlattice $F$ of continuous functions on $X$. Therefore, bypassing any subtleties that will arise when considering extensions or restrictions of order continuous linear functionals, the triviality of the order continuous duals of a particular vector lattice $E$ of continuous functions is, when obtained from the above theorem, inherited by all its superlattices $F$ of continuous functions.

We shall now combine the remainder of the discussion with an overview of the paper.

In Section~\ref{sec:preliminaries}, we introduce basic notation and establish a few auxiliary results on the role of the countable chain condition and isolated points regarding order continuous duals of vector lattices of continuous functions. The core statement of the main triviality theorem above  is that the order continuous dual of $E$ is trivial, but if $X$ satisfies the countable chain condition then the order continuous dual and the $\sigma$-order continuous dual coincide. Hence the $\sigma$-order continuous dual is then also trivial.

Section~\ref{sec:main_triviality_theorem} contains the proof of the main triviality theorem above.

It is clear from the main triviality theorem that attention should be paid to resolvable spaces, and this is the topic of Section~\ref{sec:resolvable_spaces}. It is easy to see that a non-empty resolvable topological space has no isolated points. Indeed, if $X$ is resolvable and $x\in X$ is such that $\{x\}$ is an open subset of $X$, then $\{x\}\cap D\neq \emptyset$ and $\{x\}\cap \left(X\setminus D\right)\neq \emptyset$, so that $x\in D$ as well as $x\in X\setminus D$. If $X$ avoids this obstruction by having no isolated points, then it is known to be resolvable in a number of cases: if $X$ is a metric space (see~\cite[Theorem~41]{Hewitt1943} and \cite[Theorem~3.7]{ComfortGarcia-Ferreira1994}), if $X$ is a locally compact Hausdorff space (see~\cite[Theorem~3.7]{ComfortGarcia-Ferreira1994}), if $X$ is first countable (see~\cite[Corollary to Theorem~48]{Hewitt1943}), if $X$ is a regular Hausdorff countably compact topological space (see~\cite[Theorem~6.7]{ComfortGarcia-Ferreira1994}), and, if one assumes the Axiom of Constructibility, i.e.\ if one assumes that `$V=L$', then this is also true if $X$ is a Baire topological space (see~\cite[Theorem~7.1]{ComfortGarcia-Ferreira1994}). For a useful application of the main triviality theorem, however, it is needed that $X$ has a resolving subset that can be decomposed as $\bigcup_{n\in\mathbb N}\Gamma_n$, in such a manner that the $\Gamma_n$ are `naturally' related to the topology of $X$ and to separation properties of $E$. The proofs in the papers just cited give no information about this being true or not in the  generality of the pertaining result, but Section~\ref{sec:resolvable_spaces} shows that this is indeed the case in a number of reasonably familiar situations.  The $\Gamma_n$ in that section are typically closed nowhere dense subsets of $X$, and in that case the requirement under (2) in the theorem above specialises to a (mild) separation property for points and closed nowhere dense subsets.


In Section~\ref{sec:triviality_of_order_continuous_duals}, the material from Section~\ref{sec:resolvable_spaces} on resolvable spaces and the preparatory results from Section~\ref{sec:preliminaries} are combined with the main triviality theorem from Section~\ref{sec:main_triviality_theorem}. We thus obtain results on the triviality of order continuous duals in a number of not uncommon contexts. It is our hope that the results in this section are accessible without knowledge of the rest of the paper, apart from an occasional glance at Section~\ref{sec:preliminaries} for notations and definitions.

In Section~\ref{sec:measure_theoretic_interpretation}, the triviality of the order continuous duals of two common vector lattices of continuous functions, as established in Section~\ref{sec:triviality_of_order_continuous_duals}, is interpreted in measure-theoretic terms.  The order continuous linear functionals on these two vector lattices of continuous functions correspond to so-called normal measures on certain ($\sigma$-)algebras of subsets of the underlying topological space. Hence the triviality of the order continuous duals can be interpreted as the absence of non-zero normal measures. We do not claim to give a complete overview here, and any omission in this section is entirely unintentional.


\section{Preliminaries}\label{sec:preliminaries}

In this section, we first introduce some notations and conventions. After that, we collect a few  facts about the countable chain condition for a topological space and its relation with order continuous duals of vector lattices of continuous functions, and about the role of isolated points for order continuous duals of such vector lattices. In Section~\ref{sec:triviality_of_order_continuous_duals}, these facts will be combined with the main triviality result from Section~\ref{sec:main_triviality_theorem} and the material on resolvable spaces in Section~\ref{sec:resolvable_spaces}.

\medskip

We set $\mathbb N\coloneqq\{1,2,3,\ldots\}$.

All vector spaces are over the real numbers, unless otherwise specified. If $E$ is a vector lattice, then $\orderdual{E}$ denotes its order dual. As in \cite[p.~123]{Zaanen1983RSII}, $\varphi\in\orderdual{E}$ is called \emph{$\sigma$-order continuous} if $\inf\{\,\lvert\varphi(f_n)\rvert : n\in\mathbb N\,\}=0$ whenever $(f_n)$ is a sequence in $E$ such that $f_n\downarrow 0$ in $E$, and, as already mentioned in Section~\ref{sec:introduction_and_overview}, $\varphi$ is called \emph{order continuous} if $\inf\{\,\lvert\varphi(f_i)\rvert : i\in I\,\}=0$ whenever $(f_i)_{i\in I}$ is a net in $E$ such that $f_i\downarrow 0$ in $E$. We denote the $\sigma$-order continuous linear functionals on $E$ and the order continuous linear functionals on $E$ by $\ordercontc{E}$ and $\ordercontn{E}$, respectively. Obviously, $\ordercontn{E}\subseteq\ordercontc{E}$. Furthermore, if $\varphi\in\orderdual{E}$, then, by \cite[Lemma~84.1]{Zaanen1983RSII}, $\varphi\in\ordercontc{E}$ if and only if $\varphi^+\in\ordercontc{E}$ and $\varphi^-\in\ordercontc{E}$, and $\varphi\in\ordercontn{E}$ if and only if $\varphi^+\in\ordercontn{E}$ and $\varphi^-\in\ordercontn{E}$.

A topological space will simply be called a \emph{space}. There are no implicit general suppositions concerning our spaces.  A \emph{$\text{T}_1$ space} is a space in which singletons are closed subsets. An \emph{isolated point} of a space is a point such that the corresponding singleton is an open subset. For a space $X$ and a point $x\in X$, we let $\mathcal{V}_x$ denote the collection of open neighbourhoods of $x$ in $X$. The closure of a subset $A$ of $X$ is denoted by $\overline{A}$. If $X$ is a space, then $\zerofunction$ denotes the function on $X$ that is identically zero, and $\onefunction$ denotes the function on $X$ that is identically one. If $f$ is a continuous function on a space $X$, we denote by $\zeroset f$ the zero set of $f$; i.e.\ $\zeroset f \coloneqq \{\,x\in X : f(x)=0\,\}$. When we speak of a \emph{vector lattice of continuous functions} on a space  $X$,  we shall always suppose that the partial ordering and the lattice operations are pointwise. We let $\contX$, $\contbX$, $\contoX$, and $\contcX$ denote the vector lattices of the continuous functions on $X$, the bounded continuous functions on $X$, the continuous functions on $X$ that vanish at infinity, and the compactly supported continuous functions on $X$, respectively. If $X$ is a metric space, then $\LipbX$ denotes the bounded Lipschitz functions on $X$. This is likewise a vector lattice of continuous functions; see e.g. \cite[Proposition~1.5.5] {Weaver1999LipschitzAlgebras}.

We recall that a space is a \emph{Baire space} if the intersection of an at most countably infinite collection of dense open subsets is dense as well. Equivalently, a space is a Baire space if the union of an at most countably infinite collection of closed nowhere dense subsets has empty interior. Complete metric spaces and locally compact Hausdorff spaces are Baire spaces; there exist metrizable Baire spaces which are not completely metrizable.

\medskip

We now turn to the countable chain condition and its relation with order con\-tin\-u\-ous duals of vector lattices of continuous functions.

\begin{defn}
A space $X$ is said to \emph{satisfy the countable chain condition}, or to \emph{satisfy CCC}, if every collection of non-empty pairwise disjoint open subsets of $X$ is at most countably infinite.
\end{defn}

We include the short proof of the following folklore result for the convenience of the reader.

\begin{lem}\label{lem:separable_implies_ccc}
Every separable space satisfies CCC.
\end{lem}

\begin{proof}
Suppose that $X$ is a separable space with a (possibly finite) dense subset $D=\{\,d_n : n\in\mathbb{N}\,\}$.  Consider a collection $\mathcal{C}$ of non-empty pairwise disjoint open subsets of $X$.  Since $D$ is dense in $X$, it follows that $U\cap D\neq\emptyset$ for every $U\in\mathcal{C}$.  For each $U\in\mathcal{C}$, let $n_{U}=\min\{\,n\in\mathbb{N} : d_{n}\in U\,\}$.  Since the members of $\mathcal{C}$ are pairwise disjoint, it follows that the map that sends $C\in\mathcal C$ to $n_U\in\mathbb N$
is injective.  Hence $\mathcal{C}$ is at most countably infinite.
\end{proof}

We recall that a vector lattice $E$ is \emph{order separable} if every subset of $E$ that has a supremum contains an at most countably infinite subset that has the same supremum.

For Archimedean vector lattices, there are several equivalent alternate char\-ac\-ter\-i\-sa\-tions of order separability available; see \cite[Theorem 29.3]{LuxemburgZaanen1971RSI}. One of these is the property that every subset of $E^+$ that is bounded from above and that consists of pairwise disjoint elements is at most countably infinite. It is this property that is used to establish the following result, which is a rather obvious generalisation of \cite[Theorem~10.3~(i)]{dePagterHuijsmans1980II}. We include the short proof for the sake of completeness.

\begin{lem}\label{lem:order_separability}
Let $X$ be a non-empty space that satisfies CCC. Then every vector lattice of continuous functions on $X$  is order separable.
\end{lem}

\begin{proof}
It is clearly sufficient to prove that $\contX$ is order separable. Suppose that $G\subseteq \contX^+\setminus\{0\}$ is bounded from above and that it consists of pairwise disjoint elements. It is sufficient to prove that every such $G$ is at most countably infinite. For each $v\in G$, the set $U_{v}=X\setminus \zeroset v$ is open and non-empty.  Since $G$ consists of pairwise disjoint elements, $U_{v_1}\cap U_{v_2}=\emptyset$ if $v_1,v_2\in G$ are different. This implies that $U_{v_1}\neq U_{v_2}$ if $v_1,v_2\in G$ are different. Consequently, the map that sends $v\in G$ to $U_v$ is a bijection between $G$ and $\mathcal{U}_{G}=\{\,U_v : v\in G\,\}$. Since ${\mathcal U}_G$ consists of non-empty pairwise disjoint open subsets of $X$, and since $X$ satisfies CCC, ${\mathcal U}_G$ is at most countably infinite. Hence the same holds for $G$.
%
\end{proof}

If $E$ is an order separable vector lattice, then $\ordercontn{E}\!=\ordercontc{E}$; see \cite[Theorem~84.4~(i)]{Zaanen1983RSII}. Combining this with Lemmas~\ref{lem:separable_implies_ccc} and~\ref{lem:order_separability} we have the following.

\begin{pro}\label{prop:equality_of_order_duals}
Let $X$ be a non-empty space, and let $E$ be a vector lattice of continuous functions on $X$.  If $X$ satisfies CCC \textup{(}in particular, if $X$ is separable\textup{)}, then $\ordercontn{E}=\ordercontc{E}$.
\end{pro}

We conclude with a result pointing out the special role of isolated points for order continuous duals of vector lattices of continuous functions, which the reader may readily verify upon noting that the characteristic function of an isolated point is continuous.

\begin{lem}\label{lem:evaluation_functional}
Let $X$ be a space that has an isolated point $x_0$, and suppose that $E$ is a vector lattice of continuous functions on $X$ that contains the characteristic function $\chi_{\{x_0\}}$ of $\{x_0\}$. Define the evaluation map $\delta_{x_0}:E\to\mathbb R$ by setting $\delta_{x_0}(f)\coloneqq f(x_0)$. Then  $\{0\}\neq \{\delta_{x_0}\}\subseteq\left(\ordercontn{E}\right)^+\subseteq\left(\ordercontc{E}\right)^+$.
\end{lem}


\section{Main triviality theorem}\label{sec:main_triviality_theorem}

The following result is the technical heart of the current paper. The proof is ultimately inspired by ideas in  \cite[Example~21.6 (ii)]{Zaanen1997Introduction}, where it is established that $\ordercontc{\cont([0,1])}=\{0\}$.


\begin{thm}\label{thm:main_triviality_theorem}
Let $X$ be a non-empty space, and let $E$ be a vector lattice of continuous functions on $X$.  Suppose that there exists an at most countably infinite collection $\{\,\Gamma_n : n\in\mathbb N\,\}$ of non-empty subsets of $X$ such that
 \begin{enumerate}[(1)]
    \item $\bigcup_{n\in\mathbb N}\Gamma_n$ resolves $X$, and
    \item for every $x\in X\setminus \bigcup_{n\in\mathbb N}\Gamma_n$ and every $k\in\mathbb N$ , there exists an element $u$ of $E^+$ such that $u(x)=1$ and $u(y)=0$ for all $y\in \Gamma_k$.
\end{enumerate}
Then $\ordercontn{E}=\{0\}$.

If, in addition, $X$ satisfies CCC \textup{(}in particular, if, in addition, $X$ is separable\textup{)}, then also $\ordercontc{E}=\{0\}$.
\end{thm}

%

\begin{proof}
As a preparatory reduction, we note that we may (and shall) assume that also
$\Gamma_1\subseteq\Gamma_2\subseteq\Gamma_3\subseteq\cdots $. Indeed, for $n\in\mathbb N$, set $\widetilde\Gamma_n\coloneqq\bigcup_{j=1}^n \Gamma_j$. Then $\bigcup_{n\in\mathbb N}\widetilde\Gamma_n=\bigcup_{n\in\mathbb N}\Gamma_n$ resolves $X$. If $x\in X\setminus\bigcup_{n\in\mathbb N}\widetilde\Gamma_n= X\setminus\bigcup_{n\in\mathbb N}\Gamma_n$, then, for each $j\in\mathbb N$, there exists $u_{jx}\in E^+$ such that $u_{jx}(x)=1$ and $u_{jx}(y)=0$ for all $y\in\Gamma_j$. For $k\in\mathbb N$, set $u_{kx}\coloneqq\bigwedge_{j=1}^k u_{jx}$. Then $u_{kx}\in E^+$, $u_{kx}(x)=1$, and $u_{kx}(y)=0$ for all $y\in\widetilde\Gamma_k$. Hence we may replace the $\Gamma_n$ with the  $\widetilde\Gamma_n$, and the latter form a non-decreasing chain.

After this reduction, we start with the triviality of $\ordercontn{E}$.

Arguing by contradiction, suppose that $\ordercontn{E}\neq \{0\}$. Then there exist $\varphi\in \ordercontn{E}$ and $e\in E$ such that $\varphi\geq 0$, $e\geq 0$, and $\varphi(e)=1$. Choose and fix such $\varphi$ and $e$ for the remainder of the proof.

The first step of the proof consists of introducing an auxiliary set function on the power set of $X$, as follows.

For every $n\in\mathbb N$ and every subset $A$ of $X$, set
\[
\rho_n (A) \coloneqq
\begin{cases}
1 & \textup{if } A\cap \Gamma_n\neq \emptyset;\\
0 & \textup{if } A\cap \Gamma_n=\emptyset,
\end{cases}
\]
and set
\[
\rho(A) = \sum_{n=1}^\infty \frac{1}{2^n}\rho_{n}(A).
\]
Since $\bigcup_{n\in\mathbb N}\Gamma_n$ is dense in $X$, it follows that $\rho(U)>0$ for every non-empty open subset $U$ of $X$. Furthermore, $\rho:2^X\to[0,1]$ is monotone. Both these properties are essential in the sequel of the proof; we shall encounter a similar auxiliary function on the non-empty open subsets of a space in Proposition~\ref{prop:rho_properties}.

The second step of the proof consists of constructing elements of $E$ that are indexed by $m\in\mathbb N$ and that are suitably related to the $\Gamma_n$, to $\rho$, and to $\varphi$. We shall use the order continuity of $\varphi$ here. The construction of these elements, which will be denoted by $\widetilde v_{mF_m}$ below, is as follows.

Fix $m\in \mathbb N$. We choose and fix $N_m\in\mathbb N$ such that
\begin{equation*}
\sum_{n=N_m+1}^\infty\frac{1}{2^n}\leq\frac{1}{m},
\end{equation*}
and set
\begin{equation}\label{eq:v_m_definition}
A_m=X\setminus \Gamma_{N_m}.
\end{equation}
Due to the non-decreasing nature of the chain $\{\,\Gamma_n : n\in\mathbb{N}\,\}$, we have
\begin{equation}\label{eq:upper_bound_for_measure}
\rho(A_m)\leq\sum_{n=N_m+1}^\infty\frac{1}{2^n}\leq\frac{1}{m}.
\end{equation}

Let $x$ be an arbitrary element of the non-empty set $X\setminus\bigcup_{n\in\mathbb N}\Gamma_n$.
As a con\-se\-quence of the second part of the hypotheses and the fact that $e(x)\geq 0$, there exists $u_{mx}\in E^+$ such that $u_{mx}(x)=e(x)$ and $u_{mx}(y)=0$ for all $y\in \Gamma_{N_m}$. Set $v_{mx}\coloneqq e\wedge u_{mx}$.  Then $v_{mx}\in E$, and
\begin{equation}\label{eq:first_list_of_properties}
\begin{split}
&0\leq v_{mx}\leq e,\\
&v_{mx}(x)=e(x),\\
&v_{mx}(y)=0\textup{ for all }y\in \Gamma_{N_m}.
\end{split}
\end{equation}

Denote by $\mathcal{F}$ the collection of non-empty finite subsets of $X\setminus \bigcup_{n\in\mathbb N}\Gamma_n$. Then $\mathcal F$ is non-empty, since it contains the singletons $\{x\}$ for all $x\in X\setminus \bigcup_{n\in\mathbb N}\Gamma_n$. We introduce a partial ordering on $\mathcal F$ by inclusion; it is then directed. For $F\in\mathcal F$, set $\widetilde v_{mF}\coloneqq\bigvee_{x\in F}v_{mx}\in E$.  With $m$ still fixed, the subset $\{\,\widetilde v_{mF} : F\in\mathcal{F}\,\}$ of $E$ is obviously upward directed. Furthermore, \eqref{eq:first_list_of_properties} implies that, for all $F\in\mathcal F$,
\begin{equation}\label{eq:second_list_of_properties}
\begin{split}
&0\leq \widetilde v_{mF}\leq e,\\
&\widetilde v_{mF}(x)=e(x)\textup{ for all }x\in F,\\
&\widetilde v_{mF}(y)=0\textup{ for all }y\in \Gamma_{N_m}.
\end{split}
\end{equation}

We see from \eqref{eq:second_list_of_properties} that $e$ is an upper bound of $\{\,\widetilde v_{mF} : F\in\mathcal{F}\,\}$ in $E$. We claim that actually  $\widetilde v_{mF}\uparrow e$ in $E$. To see this, let $f\in E$ be an upper bound of  $\{\,\widetilde v_{mF} : F\in\mathcal{F}\,\}$ in $E$. For every $x\in X\setminus \bigcup_{n\in\mathbb N}\Gamma_n$, taking $F=\{x\}\in\mathcal F$ in \eqref{eq:second_list_of_properties} shows that $f(x)\geq\widetilde v_{m\{x\}}(x)=e(x)$. Since $X\setminus \bigcup_{n\in\mathbb N}\Gamma_n$ is dense in $X$, and since the elements of $E$ are actually continuous (which has not been used so far), it follows that $f\geq e$. This establishes our claim that $\widetilde v_{mF}\uparrow e$ in $E$.

By the order continuity of $\varphi$, it follows that $\varphi(e-\widetilde v_{mF})\downarrow 0$.  Therefore, we can choose and fix $F_{m}\in\mathcal{F}$ such that
\begin{equation}\label{eq:phi_upper_bound}
\varphi(e-\widetilde v_{mF_m})<\frac{1}{2^{m+1}}.
\end{equation}

In the third step of the proof, we combine the $\widetilde v_{mF_m}$ as they have been found for all $m\in\mathbb N$, as follows. For each $k\in\mathbb{N}$, set
\begin{equation}\label{eq:infimum_definition}
w_{k}\coloneqq\bigwedge_{m=1}^k \widetilde v_{mF_m}.
\end{equation}
Then $(w_k)$ is decreasing and bounded below by $\zerofunction$.  We claim that actually $w_k\downarrow \zerofunction$ in $E$. To see this, suppose that $w\in E$ is such that $w\leq w_k$ for all $k\in\mathbb{N}$. Then $w^+\leq w_k$ for all $k\in\mathbb{N}$. Fix $k\in\mathbb{N}$. If $x\in X$ is such that $w^+(x)> 0$, then \eqref{eq:infimum_definition} shows that certainly $\widetilde v_{kF_k}(x)>0$. In view of \eqref{eq:second_list_of_properties}, we must then have  $x\notin \Gamma_{N_k}$, or, equivalently (see \eqref{eq:v_m_definition}), we must have $x\in A_{k}$. We conclude that $\{\,x\in X : w^+(x)>0\,\}\subseteq A_k$. Since $\rho$ is monotone, we have $\rho(\{\,x\in X : w^+(x)>0\,\})\leq\rho(A_k)$. Using \eqref{eq:upper_bound_for_measure}, we conclude that $\rho(\{\,x\in X : w^+(x)>0\,\})\leq 1/k$. Since this holds for all $k\in\mathbb N$, we see that  $\rho(\{\,x\in X : w^+(x)>0\,\})=0$. As $\rho$ is strictly positive on non-empty open subsets of $X$, it follows that $\{\,x\in X : w^+(x)>0\,\}=\emptyset$; here we use the continuity of elements of $E$ again. Hence $w^+=\zerofunction$, and then $w\leq \zerofunction$.  This establishes our claim that $w_k\downarrow \zerofunction$ in $E$.

In the fourth and final step of the proof, we use the $w_k$ and the order continuity of $\varphi$ to reach a contradiction, as follows.

Since $w_k\downarrow \zerofunction$ in $E$, we have $\varphi(w_k)\downarrow 0$. On the other hand, using that $e- \widetilde v_{mF_m}\in E^+$ for all $m\in\mathbb N$, we note that, for all $k\in\mathbb N$,
\begin{align*}
e-w_{k} & =  e - \bigwedge_{m=1}^k \widetilde v_{mF_m}\\
&=\bigvee _{m=1}^k (e- \widetilde v_{mF_m})\\
& \leq \sum_{m=1}^k (e-\widetilde v_{mF_m}).
\end{align*}
Since $\varphi\geq 0$, \eqref{eq:phi_upper_bound} therefore yields that, for all $k\in\mathbb N$,
\begin{align*}
\varphi(e-w_k)&\leq \sum_{m=1}^k \varphi(e-\widetilde v_{mF_m})\\
&< \sum_{m=1}^k \frac{1}{2^{m+1}}\\
&<\frac{1}{2}.
\end{align*}
Since $\varphi(e)=1$, it follows that
\[
\varphi(w_k)>\frac{1}{2}
\]
for all $k\in\mathbb N$. This contradicts the fact that $\varphi(w_k)\downarrow 0$, and we conclude that we must have $E^\sim_n=\{0\}$.

The second part of the statement follows from the first part and Proposition~\ref{prop:equality_of_order_duals}.
\end{proof}

\section{Resolvable spaces}\label{sec:resolvable_spaces}

In view of Theorem~\ref{thm:main_triviality_theorem},  it is clearly desirable to seek resolving subsets of resolvable spaces that can be decomposed into at most countably infinitely many components that relate naturally to continuous functions. The present section is devoted to such results.



It was already observed in Section~\ref{sec:introduction_and_overview} that a non-empty resolvable space has no isolated points, and we collect a second  elementary results for future reference.

\begin{lem}\label{lem:second_elementary_resolvable_property}
If $X$ is a resolvable non-empty $\textup{T}_1$ space, then every resolving subset of $X$ is infinite.
\end{lem}

\begin{proof}
Suppose that $X$ is a non-empty $\textup{T}_1$ space and that $D\subseteq X$ resolves $X$. If $D$ is finite, then $D$ is a closed subset of $X$, so that $X=\overline{D}=D$. But then $X=\overline{X\setminus D}=\emptyset$, which is not the case.
\end{proof}

The remainder of this section is divided into four (non-disjoint) parts, covering separable spaces, metric spaces, topological vector spaces, and locally connected Baire spaces, respectively.

\subsection{Separable spaces}

The simplest examples of resolvable spaces that have a resolving subset with a `good'  decomposition are obviously the spaces that have an at most countably infinite resolving subset. In that case,  the components of the resolving subset as in part~(1) of Theorem~\ref{thm:main_triviality_theorem} can simply be taken to be singletons. Clearly, the space is then separable, and in a number of cases this is also sufficient.

\begin{pro}\label{prop:separable_uncountable_opens_is_resolvable}
Let $X$ be a separable space. If $X$ has the property that every non-empty open subset is uncountable, then every at most countably infinite dense subset of $X$ resolves $X$.
\end{pro}

\begin{proof}
Since the resolving subsets of a space are the dense subsets with empty interior, this is clear.
\end{proof}

Proposition~\ref{prop:separable_uncountable_opens_is_resolvable} applies to a number of spaces of practical interest. For example, non-empty differentiable manifolds (which are second countable by definition) and non-empty open subsets of separable real or complex topological vector spaces have an at most countably infinite resolving subset. We shall have more to say about not necessarily separable real and complex topological vector spaces in Section~\ref{subsec:topological_topological_vector_spaces}.

We continue with another category of spaces with a countable resolving subset. It covers e.g.\ non-empty separable complete metric spaces and non-empty separable locally compact Hausdorff spaces, in both cases without isolated points.

\begin{pro}\label{prop:baire_resolvable}
Let $X$ be a non-empty separable Baire $\textup{T}_1$ space that has no isolated points. Choose an at most countably infinite dense subset $D$ of $X$. Then $D$ is actually countably infinite, and $D$ resolves $X$.
\end{pro}

\begin{proof}
Since $X$ is a $\textup{T}_1$ space that has no isolated points, the sets $X\setminus\{x\}$ are open and dense in $X$ for all $x\in X$. Since $X$ is a Baire space, $\bigcap_{x\in D}X\setminus\{x\}=X\setminus D$ is dense in $X$. Hence $D$ resolves $X$. Lemma~\ref{lem:second_elementary_resolvable_property} shows that $D$ is infinite.
\end{proof}

In the proof of Theorem~\ref{thm:separable_resolvable} we shall need the final statement of the following result, which is  \cite[Lemma 3.12]{KarlovaMykhaylyuk2013}.  The proof of \cite[Lemma 3.12]{KarlovaMykhaylyuk2013}, however, actually shows that the first statement is true, and it seems worthwhile to record this.

\begin{lem}\label{lem:separable_resolvable}
Let $X$ be a resolvable space, and suppose that the subset $D$ of $X$ is dense in $X$. Then there exists a subset of $D$ that resolves $X$. In particular, if $X$ is resolvable and separable, then there exists an at most countably infinite subset of $X$ that resolves $X$.
\end{lem}

\subsection{Metric spaces}\label{subsec:metric_spaces}

A little bit more complicated than the case where the  com\-po\-nents of a resolving subset are singletons as in Section~\ref{subsec:metric_spaces}, is the case where these components are closed nowhere dense subsets. The present section contains such results in the context of metric spaces; Section~\ref{subsec:topological_topological_vector_spaces} covers real and complex topological vector spaces, and Section~\ref{subsec:locally_connected_baire_spaces} deals with locally connected Baire spaces,

The following result applies to non-empty complete metric spaces that have no isolated points.  We recall, however, that there are metrizable Baire spaces which are not completely metrizable.

\begin{pro}\label{prop:metric_baire_is_resolvable}
Let $X$ be a non-empty metric Baire space that has no isolated points. Then there exist an at most countably infinite collection $\{\,\Gamma_n : n\in\mathbb N \,\}$ of closed nowhere dense subsets of $X$ such that $\bigcup_{n\in\mathbb N}\Gamma_n$ resolves $X$.
\end{pro}

\begin{proof}
The proof of \cite[Lemma on~p.~249]{FishelPapert1964}
shows that there exist closed nowhere dense subsets $\Gamma_1,\Gamma_2,\Gamma_3,\ldots$ of $X$ such that $\bigcup_{n\in\mathbb N}\Gamma_n$ is dense in $X$. Since $X$ is a Baire space, $\bigcup_{n\in\mathbb N}\Gamma_n$ has empty interior, i.e.\ $X\setminus \bigcup_{n\in\mathbb N}\Gamma_n$ is dense in $X$.
Hence $\bigcup_{n\in\mathbb N}\Gamma_n$ resolves $X$.
\end{proof}

Euclidean $n$-space falls within the scope of Proposition~\ref{prop:separable_uncountable_opens_is_resolvable}, so that it has a very simple suitably decomposable resolving subset. A more complicated resolving subset that is still suitably decomposable could be obtained by taking  the union of all metric spheres
\[
\sphere_r(0)=\{\,x\in \mathbb{R}^n : \|x\|=r\,\}
\]
with radii $r$ in a given countably infinite dense subset of $(0,\infty)$ and centred at $0$. The following result shows that this idea can also be used in a number of (not necessarily separable) metric spaces.

\begin{pro}\label{prop:metric_space_spheres}
Let $X$ be a locally connected metric space with metric $\dist(\cdot,\cdot)$. Suppose that $X$ contains at least two points. For $x_0\in X$, set $\ball_r(x_0)\coloneqq\{\,x\in X : \dist(x_0,x)<r\,\}$ and $\sphere_r(x_0)=\{\,x\in X : \dist(x_0,x)=r\,\}$ for $r>0$, and let $\mathrm D_{x_0}=\{\,\dist(x_0,x) : x\in X\,\}$.

Suppose that there exists a point $x_0\in X$ such that \begin{enumerate}[(1)]
    \item $\overline{\ball_r(x_0)}=\{\,x\in X:\dist(x_0,x)\leq r\,\}$ for every $r>0$, and
    \item $\mathrm D_{x_0}$ is connected in $\mathbb{R}$.
\end{enumerate}
If $\mathcal D \subseteq \mathrm D_{x_0}\setminus\{0\}$ is non-empty and dense in $\mathrm D_{x_0}$ \textup{(}such $\mathcal D$ exist\textup{)}, then $\bigcup_{r\in \mathcal D} \sphere_{r}(x_0)$ is dense in $X$, and if $\widetilde{\mathcal D}\subseteq \mathrm D_{x_0}\setminus\{0\}$ is countably infinite and dense in $\mathrm D_{x_0}$ \textup{(}such $\widetilde{\mathcal D}$ exist\textup{)}, then $\bigcup_{r\in \widetilde{\mathcal D}} \sphere_{r}(x_0)$ resolves $X$.

For all $r>0$, the subset $\sphere_r(x_0)$ of $X$ is a closed nowhere dense subset of $X$.
\end{pro}

\begin{proof}
Since $X$ contains at least two points, it follows from the second assumption that $\mathrm D_{x_0}$ is an interval of the form $[0,M]$ or $[0,M)$ for some $M>0$, or equal to $[0,\infty)$. Since $\mathrm D_{x_0}$ does not reduce to $\{0\}$, sets $\mathcal D$ and $\widetilde {\mathcal D}$ as in the statement quite obviously exist.

Clearly, $\sphere_r(x_0)$ is closed for all $r>0$. Furthermore, the first part of the hypotheses implies that, for $r>0$, $\sphere_r(x_0)= \overline{\ball_r(x_0)}\setminus\ball_r(x_0)$. Hence $\sphere_r(x_0)$ has empty interior for all $r>0$, so that it is indeed a closed nowhere dense subset of $X$.

Suppose that $\mathcal D \subseteq \mathrm D_{x_0}\setminus\{0\}$ is dense in $\mathrm D_{x_0}$. We shall show that $\bigcup_{r\in \mathcal D} S_{r}(x_0)$ is dense in $X$.

Suppose that this were not the case.  Then there exists a non-empty connected open subset $U$ of $X$ such that $U\cap \sphere_r(x_0)=\emptyset$ for every $r\in\mathcal D$.  Consider any $r\in\mathcal D$. Since $\{\,\dist(x_0,x) : x\in U \,\}$ is connected in $\mathbb R$, but does not contain $r$, we have either $\dist(x_0,x)<r$ for all $x\in U$, or $\dist(x_0,x)>r$ for all $x\in U$.

Choose and fix a point $x_1\in U$. From what we have just seen, if $x\in U$ is arbitrary, then $\dist(x_0,x)<r$ whenever $r\in\mathcal D$ is such that $\dist(x_0,x_1)<r$, and $\dist(x_0,x)>r$ whenever $r\in\mathcal D$ is such that $\dist(x_0,x_1)>r$.

If $\dist(x_0,x_1)$ is an interior point of $\mathrm D_{x_0}\setminus\{0\}$, then, by the density of $\mathcal D$ in $\mathrm D_{x_0}\setminus\{0\}$, there exists a sequence $(r_n)$ in $\mathcal D$ such that $r_n\downarrow \dist(x_0,x_1)$. From what we have just observed, it follows that $\dist(x_0,x)\leq r_n$ for all $x\in U$ and all $n\in\mathbb N$. Hence $\dist(x_0,x)\leq \dist(x_0,x_1)$ for all $x\in U$. Similarly, using a sequence $(r_n)$ in $\mathcal D$ such that $r_n\uparrow \dist(x_0,x_1)$, we have $\dist(x_0,x)\geq \dist(x_0,x_1)$ for all $x\in U$. We conclude that $U\subseteq\sphere_{\dist(x_0,x_1)}(x_0)$. This set, however, has empty interior. This contradiction shows that $\dist(x_1,x_0)$ is not an interior point of $\mathrm D_{x_0}\setminus\{0\}$.

If $\dist(x_0,x_1)=0$, then consideration of a sequence in $\mathcal D$ that decreases to 0 shows that $\dist(x_0,x)=0$ for all $x\in U$. That is, $U=\{x_0\}$. This, however, implies that $\{x_0\}$ is isolated in $\mathrm D_{x_0}$, which is not the case. This contradiction shows that $\dist(x_0,x_1)\neq 0$.

If  $\mathrm D_{x_0}$  is an interval of the form $[0,M)$ for some $M>0$, or equal to $[0,\infty)$, then all possibilities for $\dist(x_0,x_1)$ have now been exhausted, and this final contradiction concludes the proof of the density of $\bigcup_{r\in \mathcal D} S_{r}(x_0)$ in $X$ in these two cases.

If $\mathrm D_{x_0}=[0,M]$ for some $M>0$, we still have to consider the possibility that $\dist(x_0, x_1)=M$. Using a sequence in $\mathcal D$ that increases to $M$ then shows that $\dist(x_0,x)\geq M$ for all $x\in U$. Since the reverse inequality is obviously also satisfied, we find $U\subseteq\sphere_M(x_0)$, which is again impossible. This exhausts all possibilities for $\dist(x_0,x_1)$ if $\mathrm D_{x_0}=[0,M]$ for some $M>0$, and this final contradiction completes the proof of the density of $\bigcup_{r\in \mathcal D} S_{r}(x_0)$ in $X$ also in the case where $\mathrm D_{x_0}=[0,M]$ for some $M>0$.

We turn to the statement concerning $\widetilde D$. Suppose that $\widetilde{\mathcal D}\subseteq \mathrm D_{x_0}\setminus\{0\}$ is countably infinite and dense in $\mathrm D_{x_0}$. Then $\bigcup_{r\in \widetilde{\mathcal D}} S_{r}(x_0)$ is dense in $X$ by what we have just shown. On the other hand, since non-empty open subsets of $\mathrm D_{x_0}$ are uncountable, $\left(\mathrm D_{x_0}\setminus\{0\}\right)\setminus \widetilde{\mathcal D}$ is also dense in $\mathrm D_{x_0}$. Again by what we have just shown, $\bigcup_{r\in \left(\mathrm D_{x_0}\setminus\{0\}\right)\setminus \widetilde{\mathcal D}}\sphere_r(x_0)$ is dense in $X$. Since evidently $X\setminus \bigcup_{r\in \widetilde{\mathcal D}} S_{r}(x_0)=\{x_0\}\cup\bigcup_{r\in \left(\mathrm D_{x_0}\setminus\{0\}\right)\setminus \widetilde{\mathcal D}}\sphere_r(x_0)$, we see that $X\setminus \bigcup_{r\in \widetilde{\mathcal D}} S_{r}(x_0)$ is dense in $X$.

Hence $\bigcup_{r\in \widetilde{\mathcal D}} S_{r}(x_0)$ resolves $X$.

\end{proof}

\subsection{Topological vector spaces}\label{subsec:topological_topological_vector_spaces}

The idea of working with metric spheres in Proposition~\ref{prop:metric_space_spheres} can be adapted to the context of real and complex topological vector spaces, where a metric sphere with radius $r>0$ is replaced with $r(\overline U\setminus U)$ for a suitable open neighbourhood $U$ of $0$. Our final result in this vein is Theorem~\ref{thm:resolving_subset_tvs}, and we start with the necessary preparations.

\begin{lem}\label{lem:tvs_invariance}
Let $X$ be a \textup{(}not necessarily Hausdorff\textup{)} real or complex topological vector space, and let $U$ be a balanced open neighbourhood of $0$ in $X$.  Set $W\coloneqq  \bigcap_{r>0}rU$.  Then $W$ is invariant under scalar multiplication.
\end{lem}
\begin{proof}
Consider $x\in W$ and a scalar $\alpha$.  Fix $r>0$.  Since $x\in W$, it follows that
\[
x\in \tfrac{r}{1+|\alpha|}U.
\]
Then
\[
\alpha x\in \tfrac{\alpha r}{1+|\alpha|}U = r \left(\tfrac{\alpha}{1+|\alpha|}U\right).
\]
Since $U$ is balanced, we have $\tfrac{\alpha}{1+|\alpha|}U\subseteq U$.  Hence $\alpha x\in rU$.  Since $r>0$ is arbitrary, it follows that $\alpha x\in W$.
\end{proof}

\begin{lem}\label{lem:tvs_disjointness}
Let $X$ be a \textup{(}not necessarily Hausdorff\textup{)} real or complex topological vector space, and suppose that $U$ is a balanced convex open neighbourhood of $0$. If $r_1, r_2\in(0,\infty)$ are such that $r_1\neq r_2$, then $r_1\left(\overline{U}\setminus U\right)\cap r_2\left(\overline{U}\setminus U\right)=\emptyset$.
\end{lem}

\begin{proof}
Consider the Minkowski functional associated with $U$, given by
\[
\mu_U(x) = \inf\{\,t>0 : x\in t U\,\}.
\]
Since $U$ is balanced, convex, and absorbing, \cite[Theorem~1.35 (c)]{rudin_FUNCTIONAL_ANALYSIS_SECOND_EDITION:1991} shows that $\mu_U$ is a seminorm. Furthermore, since $U$ is an open subset of $X$, the proof of \cite[Theorem~1.36]{rudin_FUNCTIONAL_ANALYSIS_SECOND_EDITION:1991} shows that $\mu_U$ is continuous, and that
\begin{equation}\label{eq:reconstruction}
U=\{\,x\in X: \mu_U(x)<1\,\}.
\end{equation}
Suppose now that $r_1,r_2\in(0,\infty)$ are such that $r_1<r_2$, and suppose that $x\in r_1\left(\overline{U}\setminus U\right)\cap r_2\left(\overline{U}\setminus U\right)$.
Then there exist $x_1,x_2\in \overline{U}\setminus U$ such that $r_1 x_1 = x = r_2 x_2$. By the continuity of $\mu_U$ and \eqref{eq:reconstruction}, we have $\mu_U(x_1)\leq 1$. Since $x_2 = \frac{r_1}{r_2}x_1$, we see that $\mu_U(x_2)=\frac{r_1}{r_2} \mu_U(x_1)<1$. Then \eqref{eq:reconstruction} implies that $x_2\in U$, which is a contradiction.
\end{proof}

\begin{pro}\label{prop:tvs_density}
Let $X$ be a \textup{(}not necessarily Hausdorff\textup{)} real or complex to\-po\-log\-i\-cal vector space such that $0$ has an open neighbourhood that is not the whole space. Suppose that $U$ is a balanced open neighbourhood of $X$ such that $U-U\neq X$ \textup{(}such $U$ exist\textup{)}, and suppose that $S\subseteq(0,\infty)$ is dense in $(0,\infty)$.

Then $\bigcup_{r\in S} r\left(\overline{U}\setminus U \right)$ is dense in $X$, and, for all $r>0$, the set $r\left(\overline{U}\setminus U\right) $ is a closed nowhere dense subset of $X$.
\end{pro}

\begin{proof}
There exists an open neighbourhood $W$ of $0$ such that $W\neq X$, and then the continuity of the vector space operations implies that there exists an open neighbourhood $U$ of $0$ such that $U-U\subseteq W$. Since $U$ contains a balanced open neighbourhood of $0$ (see e.g.\ \cite[Theorem 3.1~(3.5)]{treves_TOPOLOGICAL_VECTOR_SPACES_DISTRIBUTIONS_AND_KERNELS:1967}), we may suppose that $U$ is balanced. This establishes the existence of balanced open neighbourhoods $U$ of $0$ such that $U-U\neq X$. We choose and fix one such $U$.

For each $r>0$, let
\[
\Gamma_r = \overline{r U}\setminus r U = r(\overline{U}\setminus U).
\]
Each $\Gamma_r$ is a closed nowhere dense subset of $X$.  We need to show that
\begin{equation}\label{TVS Thm Claim 1}
\bigcup_{r\in S}\Gamma_r \textup{ is dense  in } X.
\end{equation}
Suppose that (\ref{TVS Thm Claim 1}) were false.  Note that $X$, being a topological vector space, is locally connected.  Therefore, in that case, there exists a non-empty connected open subset $V$ of $X$ such that
\[
V\cap\Gamma_r=\emptyset
\]
for all $r\in S$. That is,
\[
V\subseteq (r U)\cup (X\setminus \overline{r U})
\]
for all $r\in S$. Since $V$ is connected, it follows that, for all $r\in S$,
\begin{equation}\label{TVS Thm E1}
\textup{either }  V\subseteq r U \textup { or } V\subseteq X\setminus \overline{r U}.
\end{equation}
As $U$ is an open neighbourhood of $0$ and $S$ is dense in $(0,\infty)$, there exists $r \in S$ such that  $V\cap r U\neq\emptyset$.  Hence \eqref{TVS Thm E1} shows that
\begin{equation}
\{\,r\in S : V\subseteq r U\,\}\neq\emptyset,
\end{equation}
which enables us to set
\[
R=\inf\{\,r\in S : V\subseteq r U\,\}.
\]
Consider $r_1,r_2\in(0,\infty)$ such that $r_1<r_2$.  Since $U$ is balanced, we have $r_1 U\subseteq r_2 U$. Combining this with the density of $S$ in $(0,\infty)$, we see that $V\subseteq rU$ for all $r>R$. If $R>0$ (as we shall demonstrate in a moment), then, using \eqref{TVS Thm E1} and the density of $S\cap(0,R)$ in $(0,R)$, we also see that $V\subseteq X\setminus \overline{rU}$ for all $r\in(0,R)$.

We show that $R\neq 0$.  Arguing by contradiction, suppose that $R=0$. From what we have just observed, it then follows that $V\subseteq \bigcap_{r>0}rU$, so that
\[
V-V \subseteq \bigcap_{r>0}rU - \bigcap_{r>0}rU.
\]
Since $V-V$ is a neighbourhood of $0$, we conclude that $\bigcap_{r>0}rU - \bigcap_{r>0}rU$ is absorbing.  Then Lemma~\ref{lem:tvs_invariance}  implies that $\bigcap_{r>0}rU - \bigcap_{r>0}rU  = X$, so that certainly $U-U=X$. This contradiction with our choice for $U$ shows that $R>0$.

We claim that
\begin{equation}
V\subseteq \overline{RU}.\label{TVS Thm Claim 2}
\end{equation}
To see this, consider $x\in V$.  We know that $x/r\in U$ for every $r\in(R,\infty)$. Since $x/r$ converges to $x/R$ as $r\downarrow R$, we see that $x/R\ \in \overline{U}$. Hence $x\in \overline{RU}$.

We claim that
\begin{equation}
V\subseteq X\setminus RU.\label{TVS Thm Claim 3}
\end{equation}
To see this, consider $x\in RU$, so that $x=Rx_R$ for some $x_R\in U$. We know that $V\subseteq X\setminus\overline{{rU}}\subseteq X\setminus{rU}$ for all $r\in(0,R)$. Hence $rU\subseteq X\setminus V$ for all $r\in(0,R)$, so that $rx_R\in X\setminus V$ for all $r\in(0,R)$. Since $rx_R$ converges to $R x_R=x$ as $r\uparrow R$, we see that $x\in X\setminus V$. This establishes \eqref{TVS Thm Claim 3}.

It follows from (\ref{TVS Thm Claim 2}) and (\ref{TVS Thm Claim 3}) that $V\subseteq \overline{RU}\setminus RU$, contradicting the fact that $\overline{RU}\setminus RU$ has empty interior.  Hence (\ref{TVS Thm Claim 1}) must hold.

\end{proof}

\begin{thm}\label{thm:resolving_subset_tvs}
Let $X$ be a \textup{(}not necessarily Hausdorff\textup{)} real or complex topological vector space, and let $S$ be a countably infinite dense subset of $(0,\infty)$.

If $X$ is a Baire space and $0$ has an open neighbourhood that is not the whole space, let $U$ be a balanced open neighbourhood of $0$ such that $U-U\neq X$ \textup{(}such $U$ exist\textup{)}.

If $0$ has a convex open neighbourhood that is not the whole space, let $U$ be a balanced convex open neighbourhood of $X$ such that $U-U\neq X$ \textup{(}such $U$ exist\textup{)}.

Then, in both cases,  $\bigcup_{r\in S} r\left(\overline{U}\setminus U \right)$ resolves $X$ and, for all $r>0$, the set $r\left(\overline{U}\setminus U\right) $ is a closed nowhere dense subset of $X$.
\end{thm}

\begin{proof}
We start with the case where $X$ is a Baire space and $0$ has an open  neigh\-bour\-hood  that is not the whole space. Proposition~\ref{prop:tvs_density} then shows that there exists a balanced open neighbourhood $U$ of $0$ such that $U-U\neq X$ and that, for any such $U$, $\bigcup_{r\in S} r\left(\overline{U}\setminus U \right)$ is dense in $X$. It also asserts that, for all $r>0$, the set $r\left(\overline{U}\setminus U\right) $ is a closed nowhere dense subset of $X$.
Since $X$ is a Baire space and the set $r\left(\overline{U}\setminus U\right) $ is a closed nowhere dense subset of $X$ for each $r>0$, it is immediate that $X\setminus \bigcup_{r\in S} r\left(\overline{U}\setminus U \right)$ is also dense in $X$. Hence  $\bigcup_{r\in S} r\left(\overline{U}\setminus U \right)$ resolves $X$.

We turn to the case where $0$ has a convex open neighbourhood $W$ that is not the whole space. An appeal to \cite[Theorem~1.14~(b)]{rudin_FUNCTIONAL_ANALYSIS_SECOND_EDITION:1991} shows that there exists a balanced convex open neighbourhood $U$ such that $U\subseteq \frac{1}{2}W$; although the result referred to is stated in a context where $X$ is Hausdorff, the proof makes no use of this fact. Then $U-U=U+U$ since $U$ is balanced, and $U+U\subseteq \frac{1}{2}W+\frac{1}{2}W\subseteq W$ since $W$ is convex. We have thus established the existence of a balanced convex open neighbourhood $U$ of $0$ such that $U-U\neq X$.

If $U$ is any such open neighbourhood of $0$, then Proposition~\ref{prop:tvs_density} shows again that $\bigcup_{r\in S} r\left(\overline{U}\setminus U \right)$ is dense in $X$ and that, for all $r>0$, the set $r\left(\overline{U}\setminus U\right)$ is a closed nowhere dense subset of $X$.
In addition, Lemma~\ref{lem:tvs_disjointness} implies that
\[
X\setminus \bigcup_{r\in S} r\left(\overline{U}\setminus U \right)\supseteq\bigcup_{r\in(0,\infty)\setminus S}r\left(\overline{U}\setminus U \right).
\]
Since $S$ is countably infinite, $(0,\infty)\setminus S$ is also dense in $(0,\infty)$. Proposition~\ref{prop:tvs_density} therefore also shows that $\bigcup_{r\in(0,\infty)\setminus S}r\left(\overline{U}\setminus U \right)$ is dense in $X$, and then this is certainly true for $X\setminus \bigcup_{r\in S} r\left(\overline{U}\setminus U \right)$. Hence  $\bigcup_{r\in S} r\left(\overline{U}\setminus U \right)$ resolves $X$.
\end{proof}

 \subsection{Locally connected Baire spaces}\label{subsec:locally_connected_baire_spaces}

This section serves as a preparation for Section~\ref{subsec:triviality_for_locally_connected_baire_spaces}, where we shall consider the triviality of the order continuous duals of vector lattices of continuous functions on locally connected Baire spaces.

The first main result of this section is Proposition~\ref{prop:rho_properties}, which can be viewed as defining and describing the support of an order continuous linear functional on a sufficiently rich vector lattice of continuous functions on a $T_1$ space without isolated points. Here the space does not yet need to be a Baire space.

The second main result of this section is Proposition~\ref{prop:strictly_positive_order_continuous_functional_gives_resolvability}, stating that a non-empty locally connected $T_1$ Baire space admits a resolving set of the form $\bigcup_{n\in\mathbb{N}}\Gamma_n$, where each $\Gamma_n$ is a closed nowhere dense subset of $X$, provided that $\ordercontn{E}$ contains a strictly positive order continuous linear functional for a sufficiently rich vector lattice $E$ of con\-tin\-u\-ous functions on $X$.

We start with some preparatory results.

\begin{lem}\label{lem:infimum_is_zero}
Let $X$ be a space, and let $x_0$ be a point in $X$ that is not isolated. Let $E$ be a vector lattice of continuous functions on $X$, and suppose that $S\subseteq E^+$ has the property  that, for every $x\neq x_0$ in $X$, there exists $f_x\in S$ such that $f_x(x)=0$. Then $\inf S=\zerofunction$ in $E$.
\end{lem}

\begin{proof}
Suppose that $g\in E$ is such that $g\leq f$ for all $f\in S$. Then $g\leq f_x$ for all $x\neq x_0$ in $X$, so that $g(x)\leq 0$ for all $x\neq x_0$ in $X$. If $g(x_0)>0$, then there is an open neighbourhood of $x_0$ on which $g$ is strictly positive. Since $x_0$ is not isolated, this neighbourhood contains a point $x_1$ that is different from $x_0$.  In that case, however, we already know that  $g(x_1)\leq 0$, contradicting that $g(x_1)>0$. Hence $g(x_0)\leq 0$, and we see that $g\leq\zerofunction$. Since $\zerofunction$ is a lower bound for $S$, we conclude that $\inf S=\zerofunction$ in $E$.
\end{proof}

\begin{pro}\label{prop:order_convergence_to_zero}
Let $X$ be a space, and let $x_0$ be a point in $X$ that is not isolated. Let $E$ be a vector lattice of continuous functions on $X$ with the property that, for every $x\neq x_0$ in $X$, there exists $f_x\in E^+$  such that  $f_x(x)=0$ and $f_x(x_0)=1$.

Suppose that  $(f_i)_{i\in I}\subseteq E^+$ is an order bounded net in $E$ with the property that, for every $U\in\mathcal{V}_{x_0}$, there exists $i_U\in I$ such that $X\setminus\zeroset{ f_i} \subseteq U$ for all $i\geq i_U$.

Then there exists an order bounded net $(g_j)_{j\in J}\subseteq E^+$ such that $g_j\downarrow \zerofunction$ in $E$ and with the property that, for all $j_0\in J$, there exists $i_0\in I$ such that $\zerofunction
\leq f_i\leq g_{j_0}$ for all $i\geq i_0$.

Consequently, if $\varphi\in \left(\ordercontn{E}\right)^+$, then $\inf\{\,\varphi(f_i) : i\in I\,\}=0$.
\end{pro}

\begin{proof}
Choose and fix an upped bound $u$ for $(f_i)_{i\in I}$ in $E^+$. Set
\[J=\{\,j \in E^+ :  j\leq u \textup{ and there exists }i_0\in I\textup{ such that }j\geq f_i\textup{ for all }i\geq i_0\,\}.
\]
Then $u\in J$, so that $J\neq\emptyset$. We introduce a partial order on $J$ by saying that, for $j_1,j_2\in J$,  $j_1\succeq j_2$ in $J$ if $j_1\leq j_2$ in $E$. We claim that $J$ is directed. Indeed, suppose that $j_1\in J$ and $i_1\in I$ are such that $j_1\geq f_i$ for all $i\geq i_1$, and  that $j_2\in J$ and $i_2\in I$ are such that $j_2\geq f_i$ for all $i\geq i_2$. There exists $i_3\in I$ such that $i_3\geq i_1$ and $i_3\geq i_2$, and then $j_1\wedge j_2\geq f_i$ for all $i\geq i_3$. Since also $0\leq j_1\wedge j_2\leq u$, we see that $j_1\wedge j_2\in J$. The inequalities $j_1\wedge j_2\succeq j_1$ and $j_1\wedge j_2\succeq j_2$ in $J$ then show that $J$ is directed.

Consider the net $(g_j)_{j\in J}$ in $E^+$ that is defined by $g_j\coloneqq j$ for $j\in J$. It is clear that this net is order bounded and decreasing. We claim that actually $g_j\downarrow\zerofunction$ in $E$. To see this, we shall employ Lemma~\ref{lem:infimum_is_zero}. Fix $x_1\neq x_0$. Then there exists $f_{x_1}\in E^+$ such that $f_{x_1}(x_1)=0$ and $f_{x_1}(x_0)=1$. The subset $U=\{\,x\in X : f_{x_1}(x)>\frac{1}{2}\textup{ and }u(x)<u(x_0)+1\,\}$ of $X$ is an open neighbourhood of $x_0$, so, by assumption, there exists $i_U\in I$ such that $\{\,x\in X : f_i(x)\neq 0\,\}\subseteq U$ for all $i\geq i_U$.  If $i\geq i_U$ and $x\in U$, then $2(u(x_0)+1)f_{x_1}(x)>u(x_0)+1>u(x)\geq f_i(x)$, and trivially $2(u(x_0)+1)f_{x_1}(x)\geq 0=f_i(x)$ for all $i\geq i_U$ and $x\notin U$. Hence $2(u(x_0)+1)f_{x_1}\geq f_i$ for all $i\geq i_U$. Then also $[2(u(x_0)+1)f_{x_1}]\wedge u\geq f_i$ for all $i\geq i_U$, and we conclude that  $[2(u(x_0)+1)f_{x_1}]\wedge u\in J$.

Since $g_{[2(u(x_0)+1)f_{x_1}]\wedge u}(x_1)=([2(u(x_0)+1)f_{x_1}]\wedge u)(x_1)=[2(u(x_0)+1)f_{x_1}(x_1)]\wedge u(x_1)=0\wedge u(x_1)=0$, and since $x_1$ is an arbitrary point in $X$ differing from $x_0$, Lemma~\ref{lem:infimum_is_zero} shows that $\inf \{\,g_j : j\in J\,\}=\zerofunction$ in $E$.  Hence  $g_j\downarrow\zerofunction$ in $E$, as claimed.

If $j_0\in J$ is given, then by the very definition of $J$ there exists $i_0\in I$ such that $0\leq f_i\leq j_0=g_{j_0}$ for all $i\geq i_0$.

We have thus established the existence of a net $(g_j)_{j\in J}$ with the required prop\-er\-ties.

The final statement of the proposition is clear since $\varphi(g_j)\downarrow 0$ for $\varphi\in \left(\ordercontn{E}\right)^+$.
\end{proof}

\begin{rem}\label{rem:order_convergence}
The net $(f_i)_{i\in I}$ in Proposition~\ref{prop:order_convergence_to_zero} converges in order to $\zerofunction$ in $E$ in the sense of \cite[Definition~1.18]{abramovich_aliprantis_INVITATION_TO_OPERATOR_THEORY:2002}. If we knew $E$ to be Dedekind complete, then \cite[Lemma~1.19 (b)]{abramovich_aliprantis_INVITATION_TO_OPERATOR_THEORY:2002} would imply that it is also convergent in order to $\zerofunction$ in $E$ in the (non-equivalent) sense of \cite[p.~33]{aliprantis_burkinshaw_POSITIVE_OPERATORS_SPRINGER_REPRINT:2006}. Since we do not know $E$ to be Dedekind complete, we have refrained from including any (ambiguous) statement about order convergence of $(f_i)_{i\in I}$ in Proposition~\ref{prop:order_convergence_to_zero}.
\end{rem}

In the following proposition, we introduce a function $\rho$ on the collection of non-empty subsets of $X$ that bears some resemblance to the function $\rho$ on $2^X$ in the proof of Theorem~\ref{thm:main_triviality_theorem}.

\begin{pro}\label{prop:rho_properties}
Let $X$ be a non-empty $T_1$ space without isolated points, and let $E$ be a vector lattice of continuous functions on $X$ such that, for every point $x_0\in X$ and every $U\in \mathcal{V}_{x_0}$, there exist $V\in\mathcal{V}_{x_0}$  and $f\in E$ such that $\overline{V}\subseteq U$, $\zerofunction\leq f\leq\onefunction$, $f(x)=1$ for all $x\in \overline{V}$, and $X\setminus\zeroset f\subseteq U$.

Let $\varphi\in \left(\ordercontn{E}\right)^+$ be a positive order continuous linear functional on $E$.

For each non-empty open subset $U$ of $X$, define
\begin{equation}\label{eq:rho_definition}
\rho(U) \coloneqq \sup\{\,\varphi(f) : f\in E, \, \zerofunction\leq f\leq\onefunction, \textup{ and }X\setminus \zeroset f\subseteq U\,\}.
\end{equation}
Then:
 \begin{enumerate}[(1)]
    \item $0\leq \rho(U)\leq \infty$ for each non-empty open subset $U$ of $X$, and $\rho(U)\leq\rho(V)$ for all non-empty open subsets $U$ and $V$ of $X$ such that $U\subseteq V$;
    \item if $U$ and $V$ are non-empty open subsets of $X$ such that $\rho(U)=\rho(V)=0$, then also $\rho(U\cup V)=0$;
    \item for every $x_0\in X$ and $\varepsilon>0$, there exists $U\in\mathcal{V}_{x_0}$ such that $\rho(U)<\varepsilon$;
    \item the set $S_\varphi \coloneqq \{\,x\in X : \rho(U)>0 \, \textup{for all } U\in\mathcal{V}_x\,\}$ is a closed subset of $X$, and if $\varphi\neq 0$, then it has non-empty interior,
    \item for $f\in E^+$, $\varphi(f)=0$ if and only if $f(x)=0$ for all $x\in S_\varphi$;
    \item $\varphi$ is strictly positive if and only if $S_\varphi = X$.
\end{enumerate}
\end{pro}

\begin{proof}
Let $U$ be a non-empty open subset of $X$.  Since the set in the right hand side of \eqref{eq:rho_definition} contains $\varphi(\zerofunction)=0$, the first part of the first conclusion is clear.  It is also clear that $\rho$ is monotone.

Before proceeding with the remaining conclusions, we note the following for later use.  Consider a non-empty open subset $V$ of $X$ such that $\rho(V)=0$.  We claim that then $\varphi(f)=0$ whenever $f\in E^+$ is bounded on $V$ and $X\setminus \zeroset f\subseteq V$.  Indeed, since $f$ is bounded on $V$ and zero outside $V$, there exists $M>0$ such that $\zerofunction\leq f\leq M\onefunction$, or $\zerofunction\leq f/M\leq\onefunction$ . Since $X\setminus\zeroset{f/M}=X\setminus\zeroset f\subseteq V$, this implies that $0\leq\varphi(f/M)\leq\rho(V)=0$, so that $\varphi(f)=0$.

We now turn to the second conclusion. Let $U$ and $V$ be non-empty open subsets of $X$ such that $\rho(U)=\rho(V)=0$.  Then
\begin{eqnarray}
\varphi(f)=0 \textup{ whenever } f\in E, \,\zerofunction\leq f\leq\onefunction, \textup{ and }X\setminus \zeroset f\subseteq U \label{Rho U = 0 Supp}
\end{eqnarray}
and
\begin{eqnarray}
\varphi(g)=0 \textup{ whenever } g\in E, \, \zerofunction\leq g\leq\onefunction, \textup{ and }X\setminus \zeroset g\subseteq V.  \label{Rho V = 0 Supp}
\end{eqnarray}
Set $U_0 \coloneqq U\setminus \overline{V}$ and $V_0\coloneqq V\setminus \overline{U}_0$.  Then $U_0$ and $V_0$ are disjoint open subsets of $X$. Consequently, $U_0 \cap\overline{V}_0=\emptyset$ and $V_0\cap \overline{U}_0=\emptyset$. It can happen that $U_0=\emptyset$, but $V_0$ is never empty. Indeed, if $V_0=\emptyset$, then $V\subseteq\overline{U}_0=\overline{U\setminus\overline{V}}\subseteq\overline{X\setminus\overline{V}}\subseteq\overline{X\setminus V}=X\setminus V$, which contradicts that $V$ is non-empty.

It is easy to see that $U_0\cup V_0$ and $X\setminus (U\cup V)$ are disjoint subsets of $X$. Furthermore, their union $U_0\cup V_0 \cup(X\setminus(U\cup V))$ is a dense subset of $X$.  To see this, it is sufficient to show that $U\cup V\subseteq \overline{U_0 \cup V_0}$.  For this, fix $x\in U\cup V$, and consider an arbitrary $W\in\mathcal{V}_x$. We are to show that $W\cap(U_0\cap V_0)\neq\emptyset$, and for this we may evidently suppose that $W\subseteq U\cup V$.  If $W\cap U_0=\emptyset$, then $W\subseteq \overline{V}$, so that $W\cap V\neq\emptyset$.  But in this case also $W\cap \overline{U}_0 =\emptyset$, as $X\setminus W$ is a closed subset of $X$ that contains $U_0$.  It follows that $W\cap V_0=W\cap V\neq \emptyset$.  We conclude that $x\in \overline{U_0\cup V_0}$.  Therefore, $U\cup V\subseteq \overline{U_0 \cup V_0}$, as desired, so that the subset $U_0\cup V_0 \cup(X\setminus(U\cup V))$ of $X$ is indeed dense in $X$.

After these preparations, we consider a fixed $h\in E$ such that $\zerofunction\leq h\leq\onefunction$ and $X\setminus\zeroset h \subseteq U\cup V$. We shall show that $\varphi(h)=0$, as follows. It follows from the first part of the hypotheses and the fact that $E$ is a vector lattice that, for each $y\in U_0$ (if any), there exists $f_y\in E^+$ such that $f_y(y)=h(y)$, $f_y\leq h$, and $X\setminus\zeroset{f_y}\subseteq U_0$.  Likewise, for each $z\in V_0$, there exists $g_z\in E^+$ such that $g_z(z)=h(z)$, $g_z\leq h$, and $X\setminus \zeroset{g_z} \subseteq V_0$.  Denote by $\mathcal{A}$ and $\mathcal{B}$ the collections of finite subsets of $U_0$ and $V_0$, respectively. For $\emptyset\in\mathcal{A}$, we define $f_\emptyset\coloneqq\zerofunction$; for $\emptyset\in\mathcal{B}$, we define $g_\emptyset\coloneqq\zerofunction$. For a non-empty $A=\{y_1,\dotsc,y_n\}\in \mathcal{A}$ (if any) and a non-empty $B=\{z_1,\dotsc,z_k\}\in\mathcal{B}$, we define
\[
\begin{array}{l}
f_A \coloneqq f_{y_1}\vee\cdots\vee f_{y_n},\\
g_B \coloneqq g_{z_1}\vee\cdots\vee g_{z_k}.
\end{array}
\]
Then, for all $A\in\mathcal{A}$ and $B\in\mathcal{B}$, we have
\begin{eqnarray}
f_A(y)=h(y)\text{ for all } \,y\in A\text{ (if any)}, \,\, \zerofunction \leq f_A\leq h\leq\onefunction, \textup{ and } X\setminus \zeroset{f_A}\subseteq U_0\label{fA Proporeties}
\end{eqnarray}
and
\begin{eqnarray}
g_B(z)=h(z) \text{ for all } \,z\in B\text{ (if any)}, \,\, \zerofunction \leq g_B\leq h\leq\onefunction, \textup{ and }X\setminus \zeroset{g_B}\subseteq V_0.\label{gB Proporeties}
\end{eqnarray}
Since $U_0$ and $V_0$ are disjoint, we also have, for all $A\in\mathcal{A}$ and $B\in\mathcal{B}$, that $f_A(x)+g_B(x)=h(x)$ for all $x\in A\cup B$, $\zerofunction \leq f_A+g_B\leq h$, and $X\setminus \zeroset{f_A+g_B}\subseteq U_0\cup V_0$.  Clearly, the sets $\{\,f_A : A\in\mathcal{A}\,\}$ and $\{\,g_B : B\in \mathcal{B}\,\}$ are both upward directed; therefore so is $\{\,f_A+g_B : A\in\mathcal{A},\,B\in \mathcal{B}\,\}$.  Furthermore, $f_A+g_B\leq h$ for all $A\in\mathcal{A}$ and $B\in\mathcal{B}$.  We claim that, in fact, $f_A+g_B\uparrow h$ in $E$.  Suppose that $w\in E$ is such that $f_A+g_B\leq w$ for all $A\in\mathcal{A}$ and $B\in\mathcal{B}$.  Then $w\geq f_\emptyset+g_\emptyset=\zerofunction$. In particular, for $x\in X\setminus (U\cup V)$, we have $w(x)\geq 0=h(x)$. If $x\in U_0\cup V_0$, then we can choose $A\in \mathcal{A}$ and $B\in\mathcal{B}$ such that $x\in A\cup B$, in which case $w(x) \geq f_A (x) + g_B(x) = h(x)$.  We have now established that $w(x)\geq h(x)$ for all $x\in U_0\cup V_0 \cup(X\setminus(U\cup V))$.
Since $U_0\cup V_0 \cup(X\setminus(U\cup V))$ is dense in $X$, it follows that $w\geq h$. We conclude that $f_A+g_B\uparrow h$, as claimed.

It follows from \eqref{Rho U = 0 Supp}, \eqref{Rho V = 0 Supp}, \eqref{fA Proporeties}, and \eqref{gB Proporeties} that $\varphi(f_A+g_B) = \varphi(f_A)+\varphi(g_B)=0$ for all $A\in\mathcal{A}$ and $B\in\mathcal{B}$.  Since $f_A+g_B\uparrow h$, the order continuity of $\varphi$ then implies that $\varphi(h)=0$, as we intended to show.

Since this is true for all $h\in E$ such that $\zerofunction\leq h\leq\onefunction$ and $X\setminus\zeroset h \subseteq U\cup V$, we see that $\rho(U\cup V)=0$. We have thus established the second conclusion.

We turn to the third conclusion. Fix $x_0\in X$ and $\varepsilon>0$. We are to show that there exists $U\in\mathcal{V}_{x_0}$ such that $\rho(U)<\varepsilon$.  Suppose that this were not the case. Then, for every $U\in\mathcal{V}_{x_0}$, there exists $f_U\in E$ such that $\zerofunction\leq f_U\leq\onefunction$, $X\setminus\zeroset{f_U}\subseteq U$, and $\varphi(f_U)\geq \varepsilon/2$.  We shall use Proposition~\ref{prop:order_convergence_to_zero} to show that this leads to a contradiction.

First of all, if $x\in X$ is such that $x\neq x_0$, then, since $X$ is $T_1$, $X\setminus\{x\}$ is an open neighbourhood of $x_0$. The first part of the hypotheses then implies that there exists $f_x\in E^+$ such that $f_x(x)=0$ and $f_x(x_0)=1$.

The first part of the hypotheses also implies that there exist $V_0\in \mathcal{V}_{x_0}$ and $f_0\in E$ such that $\zerofunction\leq f_0\leq\onefunction$ and $f_0(x)=1$ for all $x\in \overline{V}_0$. Let $\widetilde{\mathcal{V}}_{x_0}$ denote the collection of open neighbourhoods of $x_0$ that are contained in $V_0$, ordered by reverse inclusion, and consider the net $(f_W)_{W\in\widetilde{\mathcal{V}}_{x_0}}$ in $E^+$.

This net is bounded from above in $E$ by $f_0$.

Furthermore, if $U\in\mathcal{V}_{x_0}$ is given, then $X\setminus\zeroset{f_W}\subseteq U$ for all $W\in\widetilde{\mathcal V}_{x_0}$ such that $W\subseteq U\cap V_0\in\widetilde{\mathcal V}_{x_0}$.

 We can now apply Proposition~\ref{prop:order_convergence_to_zero} to the net $(f_W)_{W\in\widetilde{\mathcal{V}}_{x_0}}$  to conclude that $\inf\{\,\varphi(f_W) : W\in\widetilde{\mathcal{V}}_{x_0}\,\}=0$. This, however,  contradicts the fact that $\varphi(f_W)\geq \varepsilon/2$ for all $W\in\widetilde{\mathcal{V}}_{x_0}$. We have thus established the third conclusion.

We now turn to the fourth conclusion.  Suppose that $x\in X\setminus S_\varphi$.  By the definition of $S_\varphi$, there exists $V\in\mathcal{V}_x$ such that $\rho(V)=0$.  Since $V$ is an open subset of $X$, we have $V\in\mathcal{V}_y$ for all $y\in V$.  Hence $V\subseteq X\setminus S_\varphi$. We conclude that $S_\varphi$ is a closed subset of $X$.

We now show that $S_\varphi$ has non-empty interior if $\varphi\neq 0$.  Suppose, to the contrary, that the interior of $S_\varphi$ were empty. Since $\varphi>0$, there exists $u\in E^+$ such that $\varphi(u)>0$.  Choose and fix such $u$.  Since $S_\varphi$ is a closed subset of $X$ with empty interior, $X\setminus S_\varphi$ is a non-empty open and dense subset of $X$.  For each $z\in X\setminus S_\varphi$, there exists $V\in\mathcal{V}_z$ such that $\rho (V)=0$ and $u$ is bounded on $V$.  Denote by $\mathcal{K}$ the non-empty set of all pairs $(z,V)$, where $z\in X\setminus S_\varphi$ and $V$ is an open neighbourhood of $z$ such that $\rho(V)=0$ and $u$ is bounded on $V$, and by $\mathcal{L}$ the non-empty collection of non-empty finite subsets of $\mathcal{K}$. It follows from the first part of the hypotheses and the fact that $E$ is a vector lattice that, for each $(z,V)\in\mathcal{K}$, there exists a function $f_{(z,V)}\in E^+$ such that
\[
f_{(z,V)}\leq u, \, f_{(z,V)}(z) = u(z), \textup{ and } X\setminus\zeroset{f_{(z,V)}}\subseteq V.
\]
If $(z,V)\in \mathcal{K}$, then $u$ is bounded on $V$. Hence $f_{(z,V)}$ is also bounded on $V$. Since $\rho(V)=0$, the preliminary remark in the beginning of this proof shows that $\varphi(f_{(z,V)})=0$.  For each $T=\{(z_1,V_1),\dotsc,(z_k,V_k)\}\in\mathcal{L}$, set
\[
f_T \coloneqq f_{(z_1,V_1)}\vee\cdots\vee f_{(z_k,V_k)}\in E^+.
\]
Then
\[
f_{T}\leq u, \, f_{T}(z_i) = u(z_i)\text{ for } i=1,\dotsc,k, \textup{ and } X\setminus\zeroset{f_{T}}\subseteq V_1\cup\ldots\cup V_k.
\]
By the second conclusion, $\rho(V_1\cup\ldots\cup V_k)=0$.  Furthermore, $u$ and therefore $f_T$ is bounded on $V_1\cup\ldots\cup V_k$. It follows from the preliminary remark in the beginning of the proof that $\varphi(f_T)=0$ for all $T\in\mathcal{L}$.

 The set $\{\,f_T : T\in\mathcal{L}\,\}$ is upward directed, and $f_T\leq u$ for all $T\in\mathcal{L}$.  We claim that $f_T \uparrow u$ in $E$.  To see this, suppose that $v\in E$ is such that $f_T \leq v$ for all $T\in \mathcal{L}$. If $x\in X\setminus S_\varphi$, then $u(x)=f_T(x)$ for some $T\in\mathcal{L}$, which implies that $u(x)=f_T(x)\leq v(x)$. Since $X\setminus S_\varphi$ is dense in $X$, it follows that $u\leq v$. Hence $f_T \uparrow u$, as claimed.

 By the order continuity of $\varphi$, we have $\varphi(f_T)\uparrow \varphi(u)$. Hence $\varphi(u)=0$, contrary to the fact that $\varphi(u)> 0$.  This contradiction concludes the proof of the fact that $S_\varphi$ has non-empty interior.

We now turn to the fifth part of the conclusion.

Let $f\in E^+$ be such that $\varphi(f)=0$, and suppose that there exists $x_0\in S_\varphi$ such that $f(x_0)\neq 0$.  Since $f\geq\zerofunction$, it follows that there exist a real number $\varepsilon_0>0$ and $V\in\mathcal{V}_{x_0}$ such that $\varepsilon_0<f(x)$ for all $x\in V$.  Since $x_0\in S_\varphi$, we have $\rho(V)>0$. Hence there exists $u\in E^+$ such that $\zerofunction\leq u\leq\onefunction$, $X\setminus\zeroset u\subseteq V$, and $\varphi(u)>0$. Then $\varepsilon_0 u\leq f$.  Indeed, if $x\in V$, then $\varepsilon_0 u(x)\leq \varepsilon_0<f(x)$, and if $x\in X\setminus V$, then $\varepsilon_0 u(x)=0\leq f(x)$. This implies that $0<\varepsilon_0\varphi(u) = \varphi(\varepsilon_0 u)\leq \varphi(f)$, so that $\varphi(f)>0$. This is a contradiction, and we conclude that $f(x)=0$ for all $x\in S_\varphi$.

Conversely, suppose that $f\in E^+$ is such that $f(x)=0$ for all $x\in S_\varphi$.  We shall show that $\varphi(f)=0$. If $S_\varphi=X$, then $f=\zerofunction$, and all is clear. If $S_\varphi\neq X$, consider all pairs $(z,V)$, where $z\in X\setminus S_\varphi$ and $V$ is an open neighbourhood of $z$ such that $V\subseteq X\setminus S_\varphi$, $\rho(V)=0$, and $f$ is bounded on $V$.  Then there exists a function $f_{(z,V)}\in E^+$ such that $f_{(z,V)} \leq f$, $f_{(z,V)}(z) = f(z)$, and $X\setminus\zeroset{f_{(z,V)}}\subseteq V$.  We note that, since $f$ is bounded on $V$, so is $f_{(z,V)}$.  As $\rho(V)=0$, it follows that $\varphi(f_{(z,V)})=0$.  Employing a method as in the proof of the fourth conclusion above, we find, using that $f(x)=0$ for all $x\in S_\varphi$, that there exists a net $f_\alpha\uparrow f$ in $E$ such that $\varphi(f_\alpha)=0$ for each $\alpha$.  Once again, the order continuity of $\varphi$ implies that $\varphi(f)=0$.

It remains to establish the sixth and final conclusion.

 Suppose that $\varphi$ is strictly positive. Choose and fix $x_0\in X$, and let $V\in\mathcal{V}_{x_0}$ be arbitrary.  By the first part of the hypotheses, there exists a function $f_V\in E$ such that $\zerofunction\leq f_V\leq\onefunction$, $f(x_0)=1$, and $X\setminus \zeroset {f_V}\subseteq V$.  By the strict positivity of $\varphi$ and the definition of $\rho$, it follows that $\rho(V)\geq \varphi(f_V)>0$.  Hence $x_0\in S_\varphi$. We conclude that $S_\varphi = X$.

 Conversely, suppose that $S_\varphi=X$, and that $u\in E$ is such that $u>0$ and $\varphi(u)=0$. Choose and fix a point $x_0$ in $X$ such that $u(x_0)>0$. By scaling $u$, we may suppose that $u(x_0)>1$. Let $V_0=\{\,x\in X: u(x)>1\,\}$. Then $V_0\in\mathcal V_{x_0}$. If $f\in E$ is such that $\zerofunction\leq f\leq\onefunction$ and $X\setminus\zeroset{f}\subseteq V_0$, then $\zerofunction\leq f\leq u$. Since $\varphi(u)=0$, this implies that $\varphi(f)=0$. This shows that $\rho(V_0)=0$. Hence $x_0\notin S_\varphi=X$, which is a contradiction. We conclude that $\varphi$ is strictly positive.
\end{proof}

Proposition~\ref{prop:rho_properties} is used to establish the following result, which makes contact with resolvability.

\begin{pro}\label{prop:strictly_positive_order_continuous_functional_gives_resolvability}
Let $X$ be a non-empty locally connected $T_1$ Baire space without isolated points, and let $E$ be a vector lattice of continuous functions on $X$ such that, for every point $x_0\in X$ and every $U\in \mathcal{V}_{x_0}$, there exist $V\in\mathcal{V}_{x_0}$  and $f\in E$ such that $\overline{V}\subseteq U$, $\zerofunction\leq f\leq\onefunction$, $f(x)=1$ for all $x\in \overline{V}$, and $X\setminus\zeroset f\subseteq U$.

Suppose that there exists a strictly positive order continuous linear functional on $E$.

Then there exists an at most countably infinite collection $\{\,\Gamma_n : n\in \mathbb N\,\}$ of closed nowhere dense subsets of $X$ such that $\bigcup_{n\in\mathbb N}\Gamma_n$ resolves $X$.

\end{pro}

\begin{proof}
Choose and fix a strictly positive order continuous linear functional $\varphi$ on $E$. For all non-empty open subsets $U$ of $X$, let $\rho(U)$ be as defined in Proposition~\ref{prop:rho_properties}.  We observe that, due to the strict positivity of $\varphi$, $\rho$ is strictly positive on the non-empty open subsets of $X$, by the sixth conclusion of Proposition~\ref{prop:rho_properties}.

Fix $n\in\mathbb N$.  By the third conclusion of Proposition~\ref{prop:rho_properties}, there exists a non-empty open subset $U$ of $X$ such that $0<\rho(U)<\frac{1}{n}$, and  Zorn's Lemma then yields a maximal (with respect to inclusion) collection $\mathcal{F}_n$ of non-empty mutually disjoint open subsets $U$ of $X$ such that $0<\rho(U)<\frac{1}{n}$ for all $U\in\mathcal F_n$.  Let $G_n = \bigcup_{U\in {\mathcal F}_n} U$.  Clearly, $G_n$ is an open subset of $X$.

We claim that $G_n$ is dense in $X$.  To see this, consider a non-empty open subset $V$ of $X$. If $x$ is a point in $V$, then there exists an open neighbourhood $W$ of $x$ such that $0<\rho(W)<\frac{1}{n}$. Since then also $0<\rho(V\cap W)<\frac{1}{n}$, the maximality of $\mathcal F_n$ implies that $(V\cap W)\cap G_n\neq\emptyset$. Hence $G_n$ is indeed dense in $X$.

We set $\Gamma_n \coloneqq X\setminus G_n$. Then $\Gamma_n$ is a closed nowhere dense subset of $X$, and we shall now show that $\bigcup_{n\in\mathbb N}\Gamma_n$ resolves $X$.

Firstly, since $X$ is a Baire space, it follows that $X\setminus \bigcup_{n\in\mathbb{N}}\Gamma_n = \bigcap_{n\in\mathbb{N}}G_n$ is dense in $X$.

Secondly, we claim that $\bigcup_{n\in\mathbb{N}}\Gamma_n$ is dense in $X$.  Arguing by contradiction, suppose that $V$ is a non-empty open subset of $X$ such that
\[
V\subseteq X\setminus \bigcup_{n\in\mathbb{N}}\Gamma_n = \bigcap_{n\in\mathbb{N}}G_n.
\]
Since $X$ is locally connected, we may suppose that $V$ is connected.  Fix $n\in\mathbb{N}$.  We have $V\subseteq G_n = \bigcup_{U\in\mathcal{F}_n}U$.  Since $V$ is connected and the members of $\mathcal{F}_n$ are mutually disjoint open subsets of $X$, there exists $U\in\mathcal{F}_n$ such that  $V\subseteq U$.  The monotonicity of $\rho$ then implies that $\rho(V)\leq \rho(U)<\frac{1}{n}$.  This is true for all $n\in\mathbb{N}$, so that $\rho(V)=0$, contradicting the strict positivity of $\rho$ on non-empty open subsets of $X$. This contradiction shows that $\bigcup_{n\in\mathbb{N}}\Gamma_n$ is dense in $X$, as claimed.

The proof is now complete.
\end{proof}

We shall see in Proposition~\ref{prop:no_strictly_positive_order_continuous_functionals} that, actually, a strictly positive order continuous linear functional on $E$ as in Proposition~\ref{prop:strictly_positive_order_continuous_functional_gives_resolvability} does not exist.

\section{Triviality of order continuous duals}\label{sec:triviality_of_order_continuous_duals}

With the material in the Sections~\ref{sec:preliminaries},~\ref{sec:main_triviality_theorem}, and~\ref{sec:resolvable_spaces} in place, it is now easy to establish a number of results on the triviality of order continuous duals of vector lattices of continuous functions in several contexts of some practical interest.

As in Section~\ref{sec:resolvable_spaces}, we divide this section into (non-disjoint) parts,  corresponding to separable spaces, metric spaces, topological vector spaces, and locally connected Baire spaces, respectively.

\subsection{Separable spaces} For separable resolvable spaces, the required separation property for $E$ in the main triviality result Theorem~\ref{thm:main_triviality_theorem} is simply one for a collection of pairs of points.

\begin{thm}\label{thm:separable_resolvable}
Let $X$ be a non-empty separable resolvable space.  Then there exists an at most countably infinite resolving subset of $X$.  In particular, if $X$ is
\begin{enumerate}[(1)]
\item a separable space with the property that every non-empty open subset is uncountable;
\item a separable $\textup{T}_1$ Baire space that has no isolated points,
\end{enumerate}
then every at most countably infinite dense subset of $X$ is a resolving subset of $X$.

If $D$ is such an at most countably infinite resolving subset of $X$, suppose that $E$ is a vector lattice of continuous functions on $X$ with the property that, for every $x\notin D$ and every $y\in D$, there exists $f\in E$ such that $f(x)=1$ and $f(y)=0$.

Then $\ordercontn{E}=\ordercontc{E}=\{0\}$.
\end{thm}

\begin{proof}
For an arbitrary non-empty separable resolvable space, Lemma~\ref{lem:separable_resolvable} shows that there exists an at most countably infinite subset of $X$ that resolves $X$.  Propositions~\ref{prop:separable_uncountable_opens_is_resolvable} and~\ref{prop:baire_resolvable} show that the spaces under~(1) and~(2), respectively, are indeed separable resolvable spaces, and that every at most countably infinite dense subset is a resolving subset. 

Let $D=\{\,d_n : n\in\mathbb{N}\,\}$ be such a subset, where repetitions are allowed. For $n\in\mathbb N$, set $\Gamma_n\coloneqq\{d_n\}$. Then Theorem~\ref{thm:main_triviality_theorem} implies that  $\ordercontn{E}=\ordercontc{E}=\{0\}$.
\end{proof}

Theorem~\ref{thm:separable_resolvable}  is as precise is possible, in the sense that the supposed separation properties of $E$ are related to one specific decomposition of one specific subset that resolves $X$. If we strengthen these separation requirements a little, and also require that $E$ contains the char\-ac\-ter\-is\-tic functions of singletons corresponding to isolated points (if any), we obtain the equivalence in the following result.

\begin{thm}\label{thm:baire_equivalences}
Let $X$ be a non-empty separable $\textup{T}_1$ Baire space, and let $E$ be a vector lattice of continuous functions on $X$ such that
\begin{enumerate}[(1)]
\item for every two different points $x$ and $y$ in $X$, there exists $f\in E^+$ such that $f(x)=1$ and $f(y)=0$, and
\item for every isolated point $x$ of $X$, the characteristic function $\chi_{\{x\}}$ of $\{x\}$ is an element of $E$.
\end{enumerate}
Then the following are equivalent:
\begin{enumerate}
\item[(1')] $\ordercontn{E}=\{0\}$;
\item[(2')] $\ordercontc{E}=\{0\}$;
\item[(3')] $X$ has no isolated points.
\end{enumerate}
\end{thm}

\begin{proof}

In view of the separability of $X$, the equivalence of (1') and (2') follows from Proposition~\ref{prop:equality_of_order_duals}.

Since we have strengthened the separation properties of $E$ to general pairs of different points, Theorem~\ref{thm:separable_resolvable} shows that (3') implies (1') and (2').

Lemma~\ref{lem:evaluation_functional} shows that (1') implies (3'), and this completes the proof.
\end{proof}

Since separable locally compact Hausdorff spaces are separable $\textup{T}_1$ Baire spaces, Theorem~\ref{thm:baire_equivalences} applies to such spaces. In that context, $\contcX$ is a natural `minimal'  sublattice of $\contX$ that meets all the requirements in Theorem~\ref{thm:separable_resolvable}, and it seems worthwhile to record explicitly the following consequence of Theorem~\ref{thm:baire_equivalences}. Quite obviously, it applies to $\contcX$, $\contoX$, $\contbX$, and $\contX$.

\begin{cor}\label{cor:locally_compact_hausdorff}
Let $X$ be a non-empty separable locally compact Hausdorff space, and let $E$ be a vector lattice of continuous functions on $X$ that contains $\contcX$.

Then the following are equivalent:
\begin{enumerate}[(1)]
\item $\ordercontn{E}=\{0\}$;
\item $\ordercontc{E}=\{0\}$;
\item $X$ has no isolated points.
\end{enumerate}
\end{cor}

\subsection{Metric spaces}

The following result is immediate from the combination of Propositions~\ref{prop:metric_baire_is_resolvable} and~\ref{prop:metric_space_spheres} and Theorem~\ref{thm:main_triviality_theorem}.

\begin{thm}\label{thm:general_metric_baire}
Let the non-empty metric space $X$ with metric $\dist(\cdot,\cdot)$ satisfy at least one of the following:
\begin{enumerate}[(1)]
\item $X$ is a Baire space that has no isolated points;
\item $X$ is locally connected, contains at least two points, and contains  a point $x_0$ such that
\begin{enumerate}
    \item $\overline{\{\,x\in X : \dist(x_0,x)<r\,\}}=\{\,x\in X:\dist(x_0,x)\leq r\,\}$ for every $r>0$, and
    \item $\{\,\dist(x_0,x) : x\in X\,\}$ is connected in $\mathbb{R}$.
\end{enumerate}
\end{enumerate}

Then there exists an at most countably infinite collection $\{\,\Gamma_n : n\in\mathbb N \,\}$ of closed nowhere dense subsets of $X$ such that $\bigcup_{n\in\mathbb N}\Gamma_n$ resolves $X$.

If $\{\,\Gamma_n : n\in\mathbb N \,\}$ is any such collection, suppose that $E$ is a vector lattice of continuous functions on $X$  with the property that, for every $x\notin \bigcup_{n\in\mathbb N}\Gamma_n$ and every $k\in\mathbb N$, there exists $f\in E$ such that $f(x)=1$ and $f(y)=0$ for all $y\in\Gamma_k$.

Then $\ordercontn{E}=\{0\}$.

If, in addition, $X$ satisfies CCC \textup{(}in particular, if, in addition, $X$ is separable{)}, then also $\ordercontc{E}=\{0\}$.
\end{thm}

If we stipulate that $E$ has a stronger separation property, then we obtain the following analogue of Theorem~\ref{thm:baire_equivalences}. It is immediate from Theorem~\ref{thm:general_metric_baire},
Lemma~\ref{lem:evaluation_functional}, and Proposition~\ref{prop:equality_of_order_duals}.

\begin{thm}\label{thm:metric_baire_equivalences}
Let $X$ be a non-empty metric Baire space, and let $E$ be a vector lattice of continuous functions on $X$ such that
\begin{enumerate}[(1)]
\item for every point $x$ and every closed nowhere dense subset $A$ of $X$ that does not contain $x$, there exists $f\in E$ such that $f(x)=1$ and $f(y)=0$ for all $y\in A$, and
\item for every isolated point $x$ of $X$, the characteristic function $\chi_{\{x\}}$ of $\{x\}$ is an element of $E$.
\end{enumerate}
Then the following are equivalent:
\begin{enumerate}
\item[(1')] $\ordercontn{E}=\{0\}$;
\item[(2')] $X$ has no isolated points.
\end{enumerate}
If, in addition, $X$ satisfies CCC \textup{(}in particular, if, in addition, $X$ is separable\textup{)}, then \textup{(1')} and \textup{(2')} are also equivalent to:
\begin{enumerate}
\item[(3')]  $\ordercontc{E}=\{0\}.$
\end{enumerate}
\end{thm}

In the context of a separable locally compact Hausdorff space $X$, we had chosen $\contcX$ as a natural `small' vector lattice of continuous functions on $X$ that has adequate separation properties and that also contains the characteristic functions of all singletons cor\-re\-spond\-ing to isolated points. In the context of a metric space $X$, we assign this role to the vector lattice $\LipbX$ of bounded Lipschitz functions on $X$. It is a consequence of \cite[Theorem~1.5.6 (a)]{Weaver1999LipschitzAlgebras} that, for every $x\in X$ and every closed subset $A$ of $X$ that does not contain $x$, there exists $f\in\LipbX$ such that  $f(x)=1$ and $f(y)=0$ for all $y \in A$. This shows that $\LipbX$ separates points and closed nowhere dense subsets as we require it, and also that it contains the characteristic functions of all singletons corresponding to isolated points. Therefore, Theorem~\ref{thm:metric_baire_equivalences} yields the following analogue of Corollary~\ref{cor:locally_compact_hausdorff}.

\begin{cor}\label{cor:metric_baire_equivalences_lipschitz}
Let $X$ be a non-empty metric Baire space, and let $E$ be a vector lattice of continuous functions on $X$ that contains $\LipbX$. Then the following are equivalent:
\begin{enumerate}[(1)]
\item $\ordercontn{E}=\{0\}$;
\item $X$ has no isolated points.
\end{enumerate}
If, in addition, $X$ satisfies CCC \textup{(}in particular, if, in addition, $X$ is separable\textup{)}, then \textup{(1)} and \textup{(2)} are also equivalent to:
\begin{enumerate}
\item[(3)]  $\ordercontc{E}=\{0\}.$
\end{enumerate}
\end{cor}

Corollary~\ref{cor:metric_baire_equivalences_lipschitz} certainly applies to all vector sublattices of $\contX$ that contain $\contbX$. Since this is a property that is preserved under the morphisms between vector sublattices that are induced by homeomorphisms between topological spaces, the following result is rather obvious. We include it nevertheless, because it shows, for example, that for Polish spaces (where separability is a consequence of their definition), (1), (2), and (3) in Corollary~\ref{cor:polish} are equivalent. This is a point that seems worth making explicit.

\begin{cor}\label{cor:polish}
Let $X$ be a non-empty space that is homeomorphic to a metric Baire space, and let $E$ be a sublattice of $\contX$ that contains $\contbX$. Then the following are equivalent:
\begin{enumerate}[(1)]
\item $\ordercontn{E}=\{0\}$;
\item $X$ has no isolated points.
\end{enumerate}
If, in addition, $X$ satisfies CCC \textup{(}in particular, if, in addition, $X$ is separable\textup{)}, then \textup{(1)} and \textup{(2)} are also equivalent to:
\begin{enumerate}
\item[(3)]  $\ordercontc{E}=\{0\}.$
\end{enumerate}
\end{cor}

\subsection{Topological vector spaces}

Every (not necessarily Hausdorff) topological vector space $X$ is completely regular; we refer the reader to \cite[p.~21, Exercise~41] {rudin_FUNCTIONAL_ANALYSIS_SECOND_EDITION:1991} or \cite[p.~188, Corollary~17]{kelley_GENERAL_TOPOLOGY_REPRINT_OF_THE_1955_EDITION:1975} for this. Hence the following result, which is immediate from Theorems~\ref{thm:resolving_subset_tvs} and Theorem~\ref{thm:main_triviality_theorem},  actually has substance, as it applies e.g.\ to the non-zero vector lattice $\contbX$.

\begin{thm}\label{thm:tvs}

Let the \textup{(}not necessarily Hausdorff\textup{)} topological vector space $X$ satisfy at least one of the following:
\begin{enumerate}[(1)]
\item $X$ is a Baire space and $0$ has an open neighbourhood that is not the whole space;
\item $0$ has a convex open neighbourhood that is not the whole space.
\end{enumerate}
Then there exists an at most countably infinite collection $\{\,\Gamma_n : n\in \mathbb N\,\}$ of closed nowhere dense subsets of $X$ such that $\bigcup_{n\in\mathbb N}\Gamma_n$ resolves $X$.

If $\{\,\Gamma_n : n\in \mathbb N\,\}$ is any such collection,  suppose that $E$ is a vector lattice of continuous functions on $X$ with the property that, for every point $x\notin\bigcup_{n\in\mathbb N}\Gamma_n$ and every $k\in\mathbb N$, there exists $f\in E$ such that $f(x)=1$ and $f(y)=0$ for all $y\in \Gamma_k$.

Then $\ordercontn{E}=\{0\}$.

If, in addition, $X$ satisfies CCC \textup{(}in particular, if, in addition, $X$ is separable\textup{)}, then also $\ordercontc{E}=\{0\}$.
\end{thm}

\subsection{Locally connected Baire spaces}\label{subsec:triviality_for_locally_connected_baire_spaces}

In this section, we obtain a triviality result for order continuous duals of sufficiently rich vector lattices of continuous functions on locally connected $\textup{T}_1$ Baire spaces that have no isolated points.

The pertinent Theorem~\ref{thm:locally_connected_Baire} follows from the combination of  Theorem~\ref{thm:main_triviality_theorem} and the material in Section~\ref{subsec:locally_connected_baire_spaces}.

We need the following two preparatory results.

The first of these shows that, actually, a strictly positive order continuous linear functional as in Proposition~\ref{prop:strictly_positive_order_continuous_functional_gives_resolvability}  does not exist.

\begin{pro}\label{prop:no_strictly_positive_order_continuous_functionals}
Let $X$ be a non-empty locally connected $\textup{T}_1$ Baire space that has no isolated points, and let $E$ be a vector lattice of continuous functions on $X$ such that,
for every point $x_0\in X$ and every $U\in \mathcal{V}_{x_0}$, there exist $V\in\mathcal{V}_{x_0}$  and $f\in E$ such that $\overline{V}\subseteq U$, $\zerofunction\leq f\leq\onefunction$, $f(x)=1$ for all $x\in \overline{V}$, and $X\setminus\zeroset f\subseteq U$.

 Then $\ordercontn{E}$ has no strictly positive elements.
\end{pro}

\begin{proof}
Suppose that, to the contrary, $\ordercontn{E}$ has a strictly positive element. Then, according to Proposition~\ref{prop:strictly_positive_order_continuous_functional_gives_resolvability}, there exists an at most countably infinite collection $\{\,\Gamma_n : n\in \mathbb N\,\}$ of closed nowhere dense subsets of $X$ such that $\bigcup_{n\in\mathbb N}\Gamma_n$ resolves $X$.

The hypothesis implies that, for every point $x$ in $X$ and every closed subset $A$ of $X$ that does not contain $x$, there exists $f\in E^+$ such that $f(x)=1$ and $f(y)=0$ for all $y\in A$.  This, combined with Theorem~\ref{thm:main_triviality_theorem}, shows that $\ordercontn{E}=\{0\}$. This contradicts the fact that $\ordercontn{E}$ has a strictly positive element.
\end{proof}

The second preparatory result relies heavily on Proposition~\ref{prop:rho_properties}. By zooming in on the interior of the `support' of a non-zero order continuous linear functional as introduced in Proposition~\ref{prop:rho_properties}, a strictly positive order continuous linear functional on a vector lattice of continuous functions on this interior is located.

\begin{lem}\label{lem:ordercontinuous_restriction}
Let $X$ be a non-empty $\textup{T}_1$ space without isolated points, and let $E$ be a vector lattice of continuous functions on $X$ such that, for every point $x_0\in X$ and every $U\in \mathcal{V}_{x_0}$, there exist $V\in\mathcal{V}_{x_0}$  and $f\in E$ such that $\overline{V}\subseteq U$, $\zerofunction\leq f\leq\onefunction$, $f(x)=1$ for all $x\in \overline{V}$, and $X\setminus\zeroset f\subseteq U$.

Suppose that there exists a non-zero positive order continuous linear functional $\varphi$ on $E$.

Let $S_\varphi$ be as defined in Proposition~\ref{prop:rho_properties}, and let $Y$ be the interior of $S_\varphi$. According to the fourth conclusion of Proposition~\ref{prop:rho_properties}, $Y$ is not empty.  Set $F\coloneqq\{\,g|_{Y} : g\in E \textup{ and }  X\setminus\zeroset g\subseteq Y\,\}$.

Then $F$ is a vector lattice of continuous functions on $Y$, and there exists a strictly positive order continuous linear functional on $F$.
\end{lem}

\begin{proof}
It is clear that $F$ is a vector lattice of continuous functions on $Y$. We proceed to establish the existence of  a strictly positive order continuous linear functional on $F$.
Let $f\in F$, and suppose that functions $g,h\in E$ are such that $X\setminus\zeroset g \subseteq Y$, $X\setminus\zeroset h \subseteq Y$, and $g|_Y=f=h|_Y$.  Then $g(x)=h(x)$ for all $x\in X$. Hence the map $\tilde{\varphi}:F\to \mathbb{R}$, given, for $f\in F$, by
\[
\tilde{\varphi}(f) \coloneqq \varphi(g),
\]
where $g\in E$ is such that $g|_Y=f$ and $X\setminus\zeroset g \subseteq Y$, is well defined.  Clearly, $\tilde{\varphi}$ is a positive linear functional on $F$.

We shall now show that $\tilde{\varphi}$ is strictly positive.  Suppose that $f\in F^+$ is such that $\tilde{\varphi}(f)=0$.  Let $g\in E$ be such that $g|_Y=f$ and $X\setminus\zeroset g \subseteq Y$.  Then $g\geq\zerofunction$ and $\varphi(g)=\tilde{\varphi}(f)=0$.  By the fifth conclusion of Proposition~\ref{prop:rho_properties}, it follows that $g(x)=0$ for all $x\in S_\varphi$. In particular, $g(x)=0$ for all $x\in Y$. Hence $f(x)=g(x)=0$ for all $x\in Y$, so that $f=\zerofunction$ in $F$. This shows that $\tilde{\varphi}$ is strictly positive on $F$.

We claim that $\tilde{\varphi}$ is order continuous on $F$. To see this, consider a net $(f_i)_{i\in I}$ in $F$ such that $f_i \downarrow \zerofunction$ in $F$.  For each $i\in I$, let $g_i\in E$ be such that $f_i = g_i |_Y$ and $X\setminus\zeroset{g_i} \subseteq Y$.  Obviously,  $g_i\geq\zerofunction$ for all $i\in I$, and $g_i\downarrow$ in $E$. We claim that, in fact, $g_i\downarrow \zerofunction$ in $E$.  To see this, let $w\in E$ satisfy $w\leq g_i$ for all $i\in I$. Then $w^+\leq g_i^+=g_i$ for all $i\in I$. This implies that $0\leq w^+(x)\leq 0$ for all $x\in X\setminus Y$, so that $w^+(x)=0$ for all $x\in X\setminus Y$. Hence $w^+|_Y \in F$.  Furthermore, if $x\in Y$, then $w^+(x)\leq g_i (x)=f_i(x)$ for all $i\in I$, so that $w^+|_Y \leq f_i$ in $F$ for all $i\in I$. Since $f_i\downarrow\zerofunction$ in $F$, this implies that $w^+|_Y\leq \zerofunction$.  Hence $w^+(x)=0$ for all $x\in Y$. We conclude that $w^+=\zerofunction$ in $E$, so that $w\leq\zerofunction$ in $E$. This establishes our claim that $g_i\downarrow\zerofunction$ in $E$. From the definition of $\tilde{\varphi}$ and the order continuity of $\varphi$ on $E$, we then have $\tilde{\varphi}(f_i)=\varphi(g_i)\downarrow 0$.  Hence $\tilde{\varphi}$ is order continuous on $F$, as claimed.
\end{proof}

It is now just a matter of putting the pieces together to obtain the following triviality result for order continuous duals of vector lattices of continuous functions on locally connected Baire spaces.

\begin{thm}\label{thm:locally_connected_Baire}
Let $X$ be a non-empty locally connected $\textup{T}_1$ Baire space that has no isolated points, and let $E$ be a vector lattice of continuous functions on $X$ such that, for every point $x_0\in X$ and every $U\in \mathcal{V}_{x_0}$, there exist $V\in\mathcal{V}_{x_0}$ and $f\in E$ such that $\overline{V}\subseteq U$, $\zerofunction\leq f\leq\onefunction$, $f(x)=1$ for all $x\in \overline{V}$, and $X\setminus\zeroset f\subseteq U$.

Then $\ordercontn{E}=\{0\}$.

If, in addition, $X$ satisfies CCC \textup{(}in particular, if, in addition, $X$ is separable\textup{)}, then also $\ordercontc{E}=\{0\}$.

\end{thm}

\begin{proof}
Suppose that, to the contrary, $\ordercontn{E}$ has a non-zero element. Then $\ordercontn{E}$ contains a non-zero positive element. With the non-empty space $Y$ and the vector lattice $F$ of continuous functions on $Y$ as defined in Lemma~\ref{lem:ordercontinuous_restriction}, Lemma~\ref{lem:ordercontinuous_restriction} asserts that there exists a strictly positive order continuous linear functional on $F$.

We claim, however, that Proposition~\ref{prop:no_strictly_positive_order_continuous_functionals} applies to the space $Y$ and the vector lattice $F$ of continuous functions on $Y$. Assuming this for the moment, it then follows from Proposition~\ref{prop:no_strictly_positive_order_continuous_functionals} that there are no strictly positive order continuous linear functionals on $F$. This contradiction implies that we must have $\ordercontn{E}=\{0\}$ after all.

It remains to verify that Proposition~\ref{prop:no_strictly_positive_order_continuous_functionals} can be applied to the space $Y$ and the vector lattice $F$ of continuous functions on $Y$.

To start with, since $Y$ is an open subspace of a locally connected $\textup{T}_1$ Baire space $X$ that has no isolated points, $Y$ inherits these properties from $X$.

Choose and fix $x_0\in Y$ and an open neighbourhood $U$ of $x_0$ in $Y$.  Then $U$ is an open neighbourhood of $x_0$ in $X$.  Therefore, there exist an open neighbourhood $V$ of $x_0$ in $X$ and a function $f\in E$ such that $\overline{V}\subseteq U$ (the closure being taken in $X$), $\zerofunction\leq f\leq\onefunction$, $f(x)=1$ for all $x\in \overline{V}$, and $X\setminus \zeroset{f}\subseteq  U$.  Clearly, $V$ is also an open neighbourhood of $x_0$ in $Y$, and the closure of $V$ in $Y$ is just the closure $\overline{V}$ of $V$ in $X$, which is contained in $U$. Since $U\subseteq Y$, $X\setminus \zeroset{f}\subseteq Y$ so that $g\coloneqq f|_Y\in F$. Obviously, $\zerofunction\leq g\leq\onefunction$. Lastly, $g$, being the restriction of $f$ to $Y$, is such that $g(x)=f(x)=1$ for all $x\in \overline{V}$, and $Y\setminus \zeroset{g} = X\setminus\zeroset{f}\subseteq U$.

We have now verified that Proposition~\ref{prop:no_strictly_positive_order_continuous_functionals} can be applied to the space $Y$ and the vector lattice $F$ of continuous functions on $Y$. The proof is complete.
\end{proof}

It should be noted that the separation properties imposed on $E$ in Theorem~\ref{thm:locally_connected_Baire} do not imply stronger separation properties of $X$ than complete regularity.  Indeed, if $X$ is a completely regular $\textup{T}_1$ space, then every closed set $F$ in $X$ is the intersection of zero sets of continuous real-valued functions on $X$. Equivalently, the collection $\{\,X\setminus\zeroset f : f\in \contX\,\}$ is a basis for the topology of $X$.  Furthermore, if $A$ and $B$ are subsets of $X$ contained in disjoint zero sets, then there exists a continuous function $f$ on $X$ such that $\zerofunction\leq f\leq\onefunction$, $f(x)=0$ for all $x\in A$, and $f(x)=1$ for all $x\in B$.  The standard reference for these facts is \cite[Sections 1.15 and 3.2]{GillmanJerison1960}.

As particular cases of Theorem~\ref{thm:locally_connected_Baire} we have the following.

\begin{cor}
Let $X$ be a non-empty completely regular and locally connected $\textup{T}_1$ Baire space that has no isolated points.  If $E$ is a vector lattice of continuous functions on $X$ that contains $\contbX$, then $\ordercontn{E}=\{0\}$.

If, in addition, $X$ satisfies CCC \textup{(}in particular, if, in addition, $X$ is separable\textup{)}, then also $\ordercontc{E}=\{0\}$.
\end{cor}

\begin{cor}
Let $X$ be a non-empty and locally connected locally compact Hausdorff space that has no isolated points.  If $E$ is a vector lattice of continuous functions on $X$ that contains $\contcX$, then $\ordercontn{E}=\{0\}$.

If, in addition, $X$ satisfies CCC \textup{(}in particular, if, in addition, $X$ is separable\textup{)}, then also $\ordercontc{E}=\{0\}$.
\end{cor}


\section{Measure-theoretic interpretation}\label{sec:measure_theoretic_interpretation}

In this section, we discuss the connection between our results on the triviality of order continuous duals of vector lattices of continuous functions on the one hand, and normal measures on the underlying topological spaces on the other hand.  Normal measures were introduced in \cite{Dixmier1951} for Stonean spaces, and later generalised to completely regular spaces in \cite{FishelPapert1964}.  In essence, a normal measure on $X$ corresponds to an order continuous linear functional on a vector lattice of continuous functions.  We consider two cases, namely, Baire measures on completely regular spaces (see \cite{Knowles1967}), and Borel measures on locally compact Hausdorff spaces (see \cite{DalesDashiellLauStrass2016}).

Let $X$ be a completely regular space, and denote by ${\bf B}^s$ the algebra of subsets of $X$ generated by $\{\,\zeroset{f} : f\in \contbX\,\}$.  A \emph{Baire measure on $X$} is a positive totally finite and (finitely) additive function $\mu:{\bf B}^s\to \mathbb{R}$, which is inner regular in the sense that
\[
\mu(B) = \sup\{\,\mu(\zeroset f) : f\in \contbX  \textup{ and } \zeroset f \subseteq B\,\}
\]
for all  $B\in{\bf B}^s$.

There is a bijective correspondence between Baire measures $\mu$ on $X$ and positive linear functionals $\varphi$ on $\contbX$, given by
\[
\varphi(f) = \int_X f\, \mathrm{d}\mu
\]
for $f\in\contbX$; see \cite{Knowles1967} and \cite{Varadarajan1961}.  A Baire  measure $\mu$ is called \emph{normal} (see \cite[Section~5]{Knowles1967}) if
\[
\int_X f_i\, d\mu \uparrow \int f\, \mathrm{d}\mu
\]
whenever $f_i\uparrow f$ in $\contbX$.  That is to say, $\mu$ is normal if and only if the associated linear functional on $\contbX$ is order continuous.

Specialising our results on the triviality of order continuous duals  of vector lattices of continuous functions to the case of $\contbX$, we therefore obtain the following.

\begin{thm}\label{thm:Completely regular normal measures}
Let $X$ be a non-empty completely regular space.  Then there are no non-zero normal Baire measures on $X$ in each of the following cases. \begin{enumerate}[(1)]
    \item $X$ is a non-empty separable resolvable space, such as \begin{enumerate}
        \item a separable space such that every non-empty open set is uncountable, or
        \item a separable $\textup{T}_1$ Baire space that has no isolated points.
        \end{enumerate}
    \item $X$ is a metric space such that \begin{enumerate}
        \item $X$ is a Baire space that has no isolated points, or
        \item $X$ is locally connected, contains at least two points, and contains  a point $x_0$ such that
            \begin{enumerate}
                \item $\overline{\{\,x\in X : \dist(x_0,x)<r\,\}}=\{\,x\in X:\dist(x_0,x)\leq r\,\}$ for every $r>0$, and
                \item $\{\,\dist(x_0,x) : x\in X\,\}$ is connected in $\mathbb{R}$.
            \end{enumerate}
        \end{enumerate}
    \item $X$ is a \textup{(}not necessarily Hausdorff\textup{)} topological vector space such that \begin{enumerate}
        \item $X$ is a Baire space and $0$ has an open neighbourhood that is not the whole space, or
        \item $0$ has a convex open neighbourhood that is not the whole space.
        \end{enumerate}
    \item $X$ is a locally connected $\textup{T}_1$ Baire space that has no isolated points.
\end{enumerate}
\end{thm}

We next consider measures on a locally compact Hausdorff space $X$. Denote by ${\bf B}$ the $\sigma$-algebra of Borel sets in $X$. As in \cite[Section 4.1]{DalesDashiellLauStrass2016}, by a \emph{regular Borel measure on $X$} we mean a finite real-valued $\sigma$-additive measure $\mu$ on $\bf B$ such that
\[
\begin{array}{lll}
\mu(B) & = & \sup\{\,\mu(K) : K\subseteq B \textup{ and } K \textup{ is compact}\} \medskip\\
& = & \inf \{\mu(U) : B\subseteq U \textup{ and } U \textup{ is open}\,\}
\end{array}
\]
for all Borel sets $B$. We let $M(X)$ denote the regular Borel measures on $X$; it is a Banach lattice when supplied with the total variation norm. If $\mu\in M(K)$, then $\mu$ induces a linear functional $\varphi_\mu:\contoX\to\mathbb R$ by setting
\[
\varphi_\mu(f)\coloneqq\int_X f\,\mathrm{d}\mu
\]
for $f\in\contoX$. The Riesz Representation Theorem \cite[Theorem~4.1.3]{DalesDashiellLauStrass2016} states that the map sending $\mu$ to $\varphi_\mu$ is an isometric lattice isomorphism between $M(K)$ and $\contoX^\sim$.

A regular Borel measure $\mu$ on $X$ is called \emph{normal} if $\varphi_\mu(f_i)\rightarrow 0$ whenever $f_i\downarrow \zerofunction$ in $\contoX$; see \cite[Definition~4.7.1]{DalesDashiellLauStrass2016}.
This implies that $\varphi_\mu$ is an order continuous linear functional on $\contoX$. Conversely, if $\mu$ is a regular Borel measure on $X$ such that $\varphi_\mu$ is an order continuous linear functional on $\contoX$, then $(\varphi_\mu)^+=\varphi_{\mu^+}$ and $(\varphi_\mu)^-=\varphi_{\mu^-}$ are also order continuous. Since these are positive linear functionals, their order continuity implies that  $\mu^+$ and $\mu^-$ are normal, and then evidently so is $\mu$. All in all, the normal measures on $X$ are in one-to-one correspondence with the order continuous linear functionals on $X$.

Specialising our  results on the triviality of order continuous duals  of vector lattices of continuous functions to the case of $\contoX$, we therefore obtain the following.

\begin{thm}\label{thm:Locally compact normal measures}
Let $X$ be a non-empty locally compact Hausdorff space that has no isolated points.  There are no non-zero normal measures on $X$ in each of the following cases. \begin{enumerate}[(1)]
    \item $X$ is separable.
    \item $X$ is metrizable.
    \item $X$ is locally connected.
\end{enumerate}
\end{thm}

The first and third conclusions in Theorem~\ref{thm:Locally compact normal measures} are also proved in \cite[Prop\-o\-si\-tion~7.20 and  Theorem~4.7.23]{DalesDashiellLauStrass2016}.  Earlier, see \cite[Proposition 19.94]{Semadeni1971}, it had been shown that $\ordercontn{\contX}=\{0\}$ if $X$ is a separable compact Hausdorff space that has no isolated points.

It is unknown to the authors whether or not the following triviality result for $\ordercontn{\contoX}$ can be obtained using the methods developed in this paper. We need a definition to be able to state it.

\begin{defn}
Let $X$ be a completely regular space. Then $X$ is an \emph{$F$-space} when
every finitely generated ideal in the algebra $\contX$ is principal.
\end{defn}

\begin{thm}\cite[Theorem 4.7.24]{DalesDashiellLauStrass2016}\label{thm:connected fspace}
Let $X$ be a non-empty connected locally compact Hausdorff $F$-space that has no isolated points.  Then there are no non-zero normal measures on $X$, equivalently, $\ordercontn{\contoX}=\{0\}$.
\end{thm}

In view of Theorems~\ref{thm:locally_connected_Baire} and~\ref{thm:connected fspace}, it is worth noting that there exists a connected compact Hausdorff space (of weight $\mathfrak{c}$) that satisfies CCC such that $\ordercontn{\contX}=\ordercontc{\contX}\neq\{0\}$.  This result is due to Plebanek \cite{Plebanek2015}; see also \cite[Theorem 4.7.26]{DalesDashiellLauStrass2016}.

\subsection*{Acknowledgements}
The results in this paper were obtained, in part, while the second author visited Leiden University from May to June 2017. He thanks the Mathematical Institute of Leiden University for their hospitality.  This visit was supported financially by the Gottfried Wilhelm Leibniz Basic Research Institute, Johannesburg, South Africa, and the National Research Foundation of South Africa, grant number 81378. The first author thanks the Department of Mathematics and Applied Mathematics of the University of Pretoria for their hospitality during a visit in March 2018.

\bibliographystyle{amsplain}
\bibliography{Mybibliography}

\end{document}